\input amstex\documentstyle{amsppt}  
\pagewidth{12.5cm}\pageheight{19cm}\magnification\magstep1
\topmatter
\title Character sheaves on disconnected groups, VIII\endtitle
\author G. Lusztig\endauthor
\address{Department of Mathematics, M.I.T., Cambridge, MA 02139}\endaddress
\thanks{Supported in part by the National Science Foundation.}\endthanks
\endtopmatter   
\document
\define\Wb{\WW^\bul}
\define\ufs{\un{\fs}}
\define\Bpq{\Bumpeq}

\define\dw{\dot w}

\define\bcl{\bar{\cl}}

\define\uD{\un D}

\define\Up{\Upsilon}

\define\bs{\bar s}

\define\hD{\hat D}
\define\hZ{\hat Z}

\define\dsv{\dashv}
\define\Lie{\text{\rm Lie }}
\define\po{\text{\rm pos}}
\define\frl{\forall}
\define\pe{\perp}
\define\si{\sim}

\define\sqc{\sqcup}

\define\qua{\quad}

\define\bK{\bar K}

\define\bZ{\bar Z}

\define\lb{\linebreak}

\define\eSb{\endSb}

\define\op{\oplus}

\define\part{\partial}
\define\em{\emptyset}

\define\ra{\rangle}

\define\iy{\infty}
\define\m{\mapsto}
\define\do{\dots}
\define\la{\langle}
\define\bsl{\backslash}

\define\sub{\subset}    
\define\bxt{\boxtimes}
\define\T{\times}
\define\ti{\tilde}
\define\nl{\newline}
\redefine\i{^{-1}}

\define\un{\underline}

\define\ot{\otimes}
\define\bbq{\bar{\QQ}_l}

\define\Ad{\text{\rm Ad}}
\define\Hom{\text{\rm Hom}}

\define\Aut{\text{\rm Aut}}

\define\ind{\text{\rm ind}}

\define\res{\text{\rm res}}

\define\tr{\text{\rm tr}}

\define\supp{\text{\rm supp}}

\define\bst{\bigstar}
\define\he{\heartsuit}

\define\a{\alpha}
\redefine\b{\beta}

\define\g{\gamma}
\redefine\d{\delta}
\define\e{\epsilon}
\define\et{\eta}
\define\io{\iota}
\redefine\o{\omega}
\define\p{\pi}
\define\ph{\phi}

\define\r{\rho}
\define\s{\sigma}
\redefine\t{\tau}

\redefine\l{\lambda}
\define\z{\zeta}
\define\x{\xi}

\define\vt{\vartheta}

\redefine\G{\Gamma}
\redefine\D{\Delta}
\define\Om{\Omega}

\define\Ph{\Phi}
\define\Ps{\Psi}

\define\dd{\bold d}

\define\kk{\bold k}

\redefine\ss{\bold s}
\redefine\tt{\bold t}

\define\II{\bold I}

\define\NN{\bold N}

\define\QQ{\bold Q}

\define\TT{\bold T}

\define\WW{\bold W}
\define\ZZ{\bold Z}

\define\ca{\Cal A}
\define\cb{\Cal B}

\define\cd{\Cal D}

\define\ci{\Cal I}

\define\ck{\Cal K}
\define\cl{\Cal L}

\define\cp{\Cal P}

\define\ct{\Cal T}

\define\cw{\Cal W}
\define\cz{\Cal Z}
\define\cx{\Cal X}

\define\fa{\frak a}
\define\fb{\frak b}
\define\fc{\frak c}
\define\fd{\frak d}
\define\fe{\frak e}
\define\ff{\frak f}

\define\fh{\frak h}

\define\fj{\frak j}
\define\fk{\frak k}

\define\fo{\frak o}
\define\fp{\frak p}
\define\fq{\frak q}
\define\fr{\frak r}
\define\fs{\frak s}
\define\ft{\frak t}

\define\fx{\frak x}
\define\fy{\frak y}

\define\fB{\frak B}
\define\fC{\frak C}
\define\fD{\frak D}
\define\fE{\frak E}

\define\fK{\frak K}

\define\fS{\frak S}

\define\fU{\frak U}

\define\ta{\ti a}
\define\tb{\ti b}
\define\tc{\ti c}
\define\td{\ti d}
\define\te{\ti e}
\define\tf{\ti f}

\define\tl{\ti l}

\define\tit{\ti t}

\define\tA{\ti A}

\define\tI{\ti I}

\define\tK{\ti K}

\define\tM{\ti M}

\define\tT{\ti T}

\define\sh{\sharp}

\define\sps{\supset}

\define\bul{\bullet}
\define\uZ{\un Z}
\define\che{\check}

\define\bfU{\bar\fU}
\define\Rf{\text{\rm R}}
\define\tfK{\ti{\fK}}
\define\bfK{\bar{\fK}}
\define\ufK{\un{\fK}}
\define\AL{Al}
\define\BBD{BBD}
\define\CU{Cu}
\define\DE{De}
\define\KA{Ka}
\define\CS{L3}
\define\AD{L9}
\define\PCS{L10}
\define\CRG{L14}
\head Introduction\endhead
Throughout this paper, $G$ denotes a fixed, not necessarily connected, reductive algebraic
group over an algebraically closed field $\kk$. This paper is a part of a series \cite{\AD}
which attempts to develop a theory of character sheaves on $G$. 

In Section 36 we associate to any subset $J$ of the set of simple reflections an algebra
$\ufK^J$ over $\QQ(v)$ (with $v$ an indeterminate) defined using certain character sheaves on
$Z_{J,G^0}$. When $J=\em$ this is the standard Iwahori-Hecke algebra attached to $\WW$. For
general $J$ we show that this algebra shares several basic features with an Iwahori-Hecke 
algebra with unequal parameters. This opens the possibility of studying Iwahori-Hecke 
algebras with unequal parameters in the framework of the theory of perverse sheaves.

In Section 37 we prove a Mackey type formula for character sheaves on $Z_{J,D}$ where
$J\sub\II$ and $D$ is a connected component of $G$. This is essentially an identity 
involving certain induction and restriction functors analogous to one in the representation
theory of reductive groups over a finite field.

In Section 38 we study a duality operation for character sheaves on a connected component
$D$ of $G$ generalizing the case $D=G^0$ considered in \cite{\CS, III, 15}. This duality 
operator is analogous to the known duality operator for representations of a reductive 
group over a finite field, see 38.12.

In Section 39 we prove a quasi-rationality property for representations of certain
extensions of an (irreducible) Weyl group. This generalizes \cite{\CS, III, (12.9.3)} and 
is a step in the proof of a key property of character sheaves (see 39.8).

{\it Errata to Part V.} 

In 27.1 replace the diagram $N_DP/U_P@<\un a<<V_1@>a>>V_2@>a'>>D$ by 
$N_DP/U_P@<\un a<<V_1@>a'>>V_2@>a''>>D$.
                                                                                 
In the commutative diagram on p.372 replace $\hZ_1@>a>>\hZ_2$, $Z_1@>b>>Z_2$, \lb
$V_1@>c>>V_2@>a'>>D$ by $\hZ_1@>a'>>\hZ_2$, $Z_1@>b'>>Z_2$, $V_1@>c'>>V_2@>a''>>D$ 
respectively.

{\it Errata to Part VI.} 

In 28.5 replace $t_0:\Ad(d\i\dw\i)(t_0)tt_0\i$ by $t_0:t\m\Ad(d\i\dw\i)(t_0)tt_0\i$.

In 28.19(a) replace $\hZ^{\cl'}_{\em,\e_D(J),D\i}$ by $\hZ^{\cl'}_{\e_D(J),D\i}$ (twice).

{\it Errata to Part VII.} 

In 32.5 replace $>b_2$ by $b_2$.

In 32.15 replace $C_u$ by $C^u$.

In 32.18(a) replace $T_{\e'(a_{r+r'}\i}$ by $T_{\e'(a_{r+r'}\i)}$.

In 32.23 replace the first $=$ by $)=$. 

After the statement of Corollary 32.23 insert: Take $D'=D\i$ so that $\D=G^0$.

In the third line of Corollary 32.24 replace $\l'$ by $\l''$.

In the fourth line of Corollary 32.24 replace $\l''$ by $\cl''$.

In 32.26, line 13, replace $tT_x$ by $\tT_x$.

In 35.6, line 8, after $h_{i-1}\i h_i\in U^*n_iU^*$, add: $(i\in[1,r])$.

\head Contents\endhead
36. The algebra $\ufK^J$.  

37. A Mackey type formula.

38. Duality.

39. Quasi-rationality.

\head 36. The algebra $\ufK^J$\endhead  
\subhead 36.1\endsubhead
We fix a connected component $D$ of $G$. We write $\e$ instead of $\e_D:\WW@>>>\WW$ (see 
26.2). Occassionally we shall consider one (or two) other component(s), say $D'$ (or 
$D',D''$); we write $\e'$ instead of $\e_{D'}$.

For $J\sub\II$ we identify $Z_{J,D}$ (see 26.2) with 
$$\{(P,gU_P);P\in\cp_J,gU_P\in D/U_P\}$$
by $(P,P',gU_P)\m(P,gU_P)$. To any $(P,gU_P)\in Z_{J,D}$ we associate an element 
$w_{P,gU_P}\in\WW$ by the following requirements. Let $z=\po(gPg\i,P)\in{}^{\e(J)}\WW^J$; 
then

(i) $w_{P,gU_P}=w_{P_1,gU_{P_1}}$ where 
$P_1=(g\i Pg\cap P)U_P\in\cp_{J\cap\e\i(\Ad(z)J)}$, 

(ii) $w_{P,gU_P}=z$ if $\e\i(\Ad(z)J)=J$. 
\nl
Note that (i),(ii) define uniquely $w_{P,gU_P}$ by induction on $|J|$: if $|J|=0$ then 
$w_{P,gU_P}$ is given by (ii); if $|J|\ge1$ and $\e\i(\Ad(z)J)=J$ then $w_{P,gU_P}$ is
again given by (ii); if $|J|\ge1$ and $\e\i(\Ad(z)J)\ne J$ then $|J\cap\e\i(\Ad(z)J)|<|J|$
and $w_{P,gU_P}$ is determined by (i) since $w_{P_1,gU_{P_1}}$ is known from the induction
hypothesis.

From definitions we see that the map $Z_{J,D}@>>>\WW$, $(P,gU_P)\m w_{P,gU_P}$ is the
composition of $\b':Z_{J,D}@>>>\ct(J,\e)$ in 26.2 (see also \cite{\PCS, 3.11}) with the 
bijection $\ct(J,\e)@>>>{}^{\e(J)}\WW$ given by \cite{\PCS, 2.4, 2.5} and with the 
inclusion ${}^{\e(J)}\WW@>>>\WW$. In particular, we have $w_{P,gU_P}\in{}^{\e(J)}\WW$ and
$\po(gPg\i,P)=\min(\WW_{\e(J)}w_{P,gU_P}\WW_J)$ for any $(P,gU_P)\in Z_{J,D}$. (For any 
$\WW_{\e(J)},\WW_J$ double coset $\Om$ in $\WW$ we denote by $\min(\Om)$ the unique element
of minimal length in $\Om$.)

For any $J\sub\II$, $w\in{}^{\e(J)}\WW$, we set 
$${}^wZ_{J,D}=\{(P,gU_P)\in Z_{J,D};w_{P,gU_P}=w\}.$$
Then $Z_{J,D}=\cup_{w\in{}^{\e(J)}\WW}{}^wZ_{J,D}$ is a partition; it is the same as the
partition 

$Z_{J,D}=\cup_{\tt\in\ct(J,\e)}{}^\tt Z_{J,D}$
\nl
in 26.2 (see also \cite{\PCS, 3.11}). In particular, each ${}^wZ_{J,D}$ is a locally 
closed, smooth subvariety of $Z_{J,D}$ stable under the $G^0$-action 
$h:(P,gU_P)\m(hPh\i,hgh\i U_{hPh\i})$.

From definitions, for $J\sub\II$, $w\in{}^{\e(J)}\WW$ we have a map 
$$\vt_{J,w}:{}^wZ_{J,D}@>>>{}^wZ_{J_1,D},(P,gU_P)\m(P_1,gU_{P_1})$$
with $J_1=J\cap\e\i(\Ad(w_0)J)$, $w_0=\min(\WW_{\e(J)}w\WW_J)$, $P_1=(g\i Pg\cap P)U_P$. 
This map may be identified with a map $\vt$ as in \cite{\PCS, 3.11}; in particular it is an
affine space bundle (see \cite{\PCS, 3.12(b)}).

\subhead 36.2\endsubhead
If $J\sub\II$ and $w\in{}^{\e(J)}\WW$ satisfies $\e\i(\Ad(w)J)=J$ (hence 
$w\in{}^{\e(J)}\WW^J$), we have ${}^wZ_{J,D}=\{(P,gU_P)\in Z_{J,D},\po(gPg\i,P)=w\}$. In
this case we pick $P\in\cp_J,P'\in\cp_{\e(J)}$ such that $\po(P',P)=w$ and a common Levi
$L$ of $P',P$. Let $\dd=\{g\in D;gLg\i=L,gPg\i=P'\}$, a connected component of the 
reductive group $N_GL$ with identity component $L$. We have a diagram
$$\dd@<pr_2<<G^0\T\dd@>j>>{}^wZ_{J,D}\tag a$$
where $j(h,g)=(hPh\i,hgh\i U_{hPh\i})$ is a principal $(P\cap P')$-bundle for a 
$(P\cap P')$-action on $G^0\T\dd$ compatible under $pr_2$ with the $(P\cap P')$-action on
$\dd$ via its quotient $(P\cap P')/U_{P\cap P'}=L$ (by conjugation in $N_GL$). If $X$ is a
character sheaf on $\dd$ then $pr_2^\bst X$ is therefore a $(P\cap P')$-equivariant simple
perverse sheaf on $G^0\T\dd$ which must be of the form $j^\bst X'$ for a well defined 
simple perverse sheaf $X'$ on ${}^wZ_{J,D}$.

The collection of simple perverse sheaves on ${}^wZ_{J,D}$ of the form $X'$ with $X$ as
above is denoted by ${}^w\hZ_{J,D}$. This collection is independent of the choice of
$P,P',L$.  Note that $X\m X'$ defines a bijection between the set of isomorphism classes of
character sheaves on $\dd$ and the set of isomorphism classes of objects in 
${}^w\hZ_{J,D}$.
 
\subhead 36.3\endsubhead
More generally, for $J\sub\II$ and $w\in{}^{\e(J)}\WW$, we define by induction on $|J|$ a
collection of simple $G^0$-equivariant perverse sheaves ${}^w\hZ_{J,D}$ on ${}^wZ_{J,D}$.

If $|J|=0$ then ${}^w\hZ_{J,D}$ is defined as in 36.2. If $|J|\ge1$ and $\e\i(\Ad(w)J)=J$,
then ${}^w\hZ_{J,D}$ is again defined as in 36.2. If $|J|\ge1$ and $\e\i(\Ad(w)J)\ne J$ 
then $\e\i(\Ad(w_0)J)\ne J$ where $w_0=\min(\WW_{\e(J)}w\WW_J)$. Thus, if 
$J_1=J\cap\e\i(\Ad(w_0)J)$ then $|J_1|<|J|$ and the class of perverse sheaves 
${}^w\hZ_{J_1,D}$ on ${}^wZ_{J_1,D}$ is defined from the induction hypothesis. By
definition, ${}^w\hZ_{J,D}$ consists of the simple perverse sheaves on ${}^wZ_{J,D}$ of the
form $\vt_{J,w}^\bst(X)$ for some $X\in{}^wZ_{J_1,D}$ (with $\vt_{J,w}$ as in 36.1). This 
completes the inductive definition of ${}^w\hZ_{J,D}$. The objects of ${}^w\hZ_{J,D}$ are 
said to be {\it character sheaves on ${}^wZ_{J,D}$.} Let $\ci_{J,w,D}$  be a set of 
representatives for the isomorphism classes of character sheaves on ${}^wZ_{J,D}$.

For $J\sub\II$ and $w\in{}^{\e(J)}\WW$, let $\cd^{cs}({}^wZ_{J,D})$ be the subcategory of
$\cd({}^wZ_{J,D})$ whose objects are those $K\in\cd({}^wZ_{J,D})$ such that for any $j$,
any simple subquotient of ${}^pH^jK$ is in ${}^w\hZ_{J,D}$. Let 
$i_{J,w}:{}^wZ_{J,D}@>>>Z_{J,D}$ be the inclusion.

(a) {\it If $K\in\cd^{cs}(Z_{J,D})$ then $i^*_{J,w}(K)\in\cd^{cs}({}^wZ_{J,D})$.}
\nl
It is enough to prove (a) for $K\in\hZ_{J,D}$. In this case (a) follows from 
\cite{\PCS, 4.12}.

For $J\sub\II$, $w\in{}^{\e(J)}\WW$ and $K\in{}^w\hZ_{J,D}$, let $K^\sh$ be the unique
simple perverse sheaf on $Z_{J,D}$ such that $i^*_{J,w}(K^\sh)=K$ and $\supp(K^\sh)$ is the
closure in $Z_{J,D}$ of $\supp(K)$.

(b) {\it We have $K^\sh\in\hZ_{J,D}$.}
\nl
This is proved in \cite{\PCS, 4.17} based on the following statement:

(c) {\it Let $Y$ be a locally closed subvariety of an algebraic variety $Y'$ and let
$i:Y@>>>Y'$ be the inclusion. Let $C\in\cd(Y)$ and let $A$ be a simple perverse sheaf on
$Y$ such that $A\dsv C$. Let $A^\sh$ be the unique simple perverse sheaf on $Y$ such that
$i^*(A^\sh)=A$ and $\supp(A^\sh)$ is the closure in $Y'$ of $\supp(A)$. Then 
$A^\sh\dsv i_!C$.}
\nl
(In {\it loc.cit.} this is applied with $Y={}^wZ_{J,D},Y'=Z_{J,D},A=K$ and an 
$C\in\cd^{cs}(Z_{J,D})$.) Since the proof of (c) is omitted in {\it loc.cit.} we give a
proof here. Let $Y_1$ be the closure of $Y$ in $Y'$. Let $Y@>i_1>>Y_1@>i_2>>Y'$ be the
inclusions. Clearly, if (c) holds for $i_1$ and $i_2$ instead of $i$ then it holds for $i$.
Thus we may assume that $i=i_1$ or $i=i_2$ that is, that $Y$ is open or closed in $Y'$. 
Assume first that $Y$ is closed in $Y'$. Then $A^\sh=i_!A$ and $i_!$ commutes with taking 
${}^pH^j$; hence $i_!A\dsv i_!C$ as desired. Next, assume that $Y$ is open in $Y'$. Then
$i^*$ commutes with taking ${}^pH^j$ and $i^*i_!C=C$. Hence for any $j$ we have 
${}^pH^j(C)=i^*({}^pH^j(i_!C))$ and for some $j$ we have $A\dsv i^*({}^pH^j(i_!C))$. Let

$0=F_0\sub F_1\sub\do\sub F_m={}^pH^j(i_!C)$
\nl
be a sequence of perverse subsheaves of ${}^pH^j(i_!C)$ such that $F_k/F_{k-1}$ is simple 
for $k\in[1,m]$. Then

$0=i^*F_0\sub i^*F_1\sub\do\sub i^*F_m=i^*({}^pH^j(i_!C))$
\nl
is a sequence of perverse subsheaves of $i^*({}^pH^j(i_!C))$ such that 
$i^*F_k/i^*F_{k-1}=i^*(F_k/F_{k-1})$ is simple or $0$ for $k\in[1,m]$. Hence 
$A\cong i^*F_k/i^*F_{k-1}=i^*(F_k/F_{k-1})$ for some $k\in[1,m]$. Then $A'=F_k/F_{k-1}$ is
a simple perverse sheaf on $Y'$ such that $A'\dsv i_!C$, $i^*A'\cong A$. We must have 
$A'\cong A^\sh$ and (c) is proved.

(d) {\it Let $K'\in\hZ_{J,D}$. There exists a unique $w\in{}^{\e(J)}\WW$ and a unique
$K\in{}^w\hZ_{J,D}$ (up to isomorphism) such that $K'\cong K^\sh$.}
\nl
This is proved in \cite{\PCS, 4.13}; it is an immediate consequence of (a).

Let $\ci_{J,D}=\sqc_{w\in{}^{\e(J)}\WW}\{K^\sh;K\in\ci_{J,w,D}\}$. From (b), (d) we see
that

(e) {\it $\ci_{J,D}$ is a set of representatives for the isomorphism classes of character
sheaves on $Z_{J,D}$.}

From (a),(b) we deduce:

(f) {\it If $w\in{}^{\e(J)}\WW$ and $K\in\cd^{cs}({}^wZ_{J,D})$, then 
$(i_{J,w})_!K\in\cd^{cs}(Z_{J,D})$.}

\subhead 36.4\endsubhead
For $J\sub J'\sub\II$ and $P\in\cp_J$ let $Q_{J',P}$ be the unique parabolic in $\cp_{J'}$
that contains $P$. We have a diagram
$$Z_{J,D}@<\fc<<Z_{J,J',D}@>\fd>>Z_{J',D}$$
where

$Z_{J,J',D}=\{(P,gU_Q);P\in\cp_J,Q=Q_{J',P},gU_Q\in D/U_Q\},$

$\fc(P,gU_Q)=(P,gU_P),\fd(P,gU_Q)=(Q,gU_Q)$. 
\nl
Define functors
$$\ff_{J,J'}:\cd(Z_{J,D})@>>>\cd(Z_{J',D}),\fe_{J,J'}:\cd(Z_{J',D})@>>>\cd(Z_{J,D})$$
by $\ff_{J,J'}A=\fd_!\fc^*A$, $\fe_{J,J'}A'=\fc_!\fd^*A'$. Now $\fd$ is proper and $\fc$ is
an affine space bundle with fibres of dimension $a=\dim\cp_J-\dim\cp_{J'}$. Hence
$\ff_{J,J'}$ commutes with Verdier duality up to a shift and a twist:

(a) $\fD\ff_{J,J'}=\ff_{J,J'}[[a]]\fD:\cd(Z_{J,D})@>>>\cd(Z_{J',D})$.
\nl
For $J\sub J'\sub J''\sub\II$ we have (see \cite{\PCS, 6.2}):

(b) $\ff_{J,J''}=\ff_{J',J''}\ff_{J,J'}$, $\fe_{J,J''}=\fe_{J,J'}\fe_{J',J''}$.
\nl
Clearly, $\ff_{J,J}=1,\fe_{J,J}=1$. 

For $J\sub J'\sub\II$ we show that the convolution bifunctor 

$\cd(Z_{J,D})\T\cd(Z_{\e(J),D'})@>>>\cd(Z_{J,D'D})$, $A,B\m A*B$,
\nl
(see 32.5) is compatible with the functors $\fe_{J,J'}$ in the following sense: for 
$A\in\cd(Z_{J,D})$, $A'\in\cd(Z_{\e(J),D'})$, we have
$$(\fe_{J,J'}A)*(\fe_{J,J'}A')=\fe_{J,J'}(A*A').\tag c$$
We have a commutative diagram
$$\CD
Z_{J',D}\T Z_{\e(J'),D'}@<b'_1<<Z'_0@>b'_2>>Z_{J',D'D}\\
@A\fd\T\fd'AA                   @Af_1AA            @A\fd''AA\\
Z_{J,J',D}\T Z_{\e(J),\e(J'),D'}@<f_2<<\uZ@>f_3>>Z_{J,J',D'D}\\
@V\fc\T\fc'VV                           @Vf_4VV     @V\fc''VV\\
Z_{J,D}\T Z_{\e(J),D'}@<b_1<<Z_0@>b_2>>Z_{J,D'D}\\    \endCD$$
Here $\fc,\fd$ are as above,

$Z_{\e(J'),D'}@<\fc'<<Z_{\e(J),\e(J'),D'}@>\fd'>>Z_{\e(J'),D'}$,
$Z_{J,D'D}@<\fc''<<Z_{J,J',D'D}@>\fd''>>Z_{J',D'D}$
\nl
are the analogous maps when $(J,J',D)$ is replaced by $(\e(J),\e(J'),D')$ or by 
$(J,J',D'D)$, $b_1,b_2$ are as in 32.5, $b'_1,b'_2$ are the analogous maps with $J$ 
replaced by $J'$,
$$\align\uZ=\{((X,P,gU_P,g'U_{P'});&X\in\cp_J,P=Q_{J',X},gU_P\in D/U_P,\\&P'=gPg\i,
g'U_{P'}\in D'/U_{P'}\},\endalign$$
$f_1(X,gU_P,g'U_{P'})=(P,gPg\i,g'gPg\i g'{}\i,gU_P,g'U_{gPg\i})$,

$f_2(X,gU_P,g'U_{P'})=((X,gU_P),(gXg\i,g'U_{gPg\i}))$,

$f_3(X,gU_P,g'U_{P'})=(X,g'gU_P)$,

$f_4(X,gU_P,g'U_{P'})=(X,gXg\i,g'gXg\i g'{}\i,gU_X,g'U_{gXg\i})$.
\nl
It is enough to show that 
$$b_{2!}b_1^*(\fc\T\fc')_!(\fd\T\fd')^*(A\bxt A')
=\fc''_!\fd''{}^*b'_{2!}b'_1{}^*(A\bxt A').$$
Since the upper right and lower left squares are cartesian we have 
$$b_1^*(\fc\T\fc')_!=f_{4!}f_2^*,\fd''{}^*b'_{2!}=f_{3!}f_1^*.$$
Hence it is enough to show that
$$b_{2!}f_{4!}f_2^*(\fd\T\fd')^*(A\bxt A')=\fc''_!f_{3!}f_1^*b'_1{}^*(A\bxt A').$$
This follows from the commutativity of the diagram above: we have 
$b_{2!}f_{4!}=\fc''_!f_{3!}$ and $f_2^*(\fd\T\fd')^*=f_1^*b'_1{}^*$. This proves (c).

If $A\in\cd(Z_{J,D})$, $B\in\cd(Z_{\e(J),D'}),C\in\cd(Z_{\e'\e(J),D''})$ then, from  
definitions, we have the associativity property:
$$(A*B)*C=A*(B*C).\tag d$$
Consider the isomorphism $\part:Z_{J,D}@>\si>>Z_{\e(J),D\i}$ as in 28.19. Note that the 
composition $Z_{J,D}@>\part>>Z_{\e(J),D\i}@>\part>>Z_{J,D}$ is the identity map.
From definitions, we have for $A,B$ as above:
$$\part_!(A*B)=(\part_!B)*(\part_!A).\tag e$$
For $J\sub\II$ we define a functor $\t:\cd(Z_{J,G^0})@>>>\cd(\text{point})$ by
$\t(C)=p'_!i^*(C)$ where $i:\cp_J@>>>Z_{J,G^0}$ is the imbedding $P\m(P,U_P)$ and 
$p':\cp_J@>>>\text{point}$ is the obvious map. We define a bifunctor 
$(:):\cd(Z_{J,D})\T\cd(Z_{J,D})@>>>\cd(\text{point})$ by $(A:B)=p_!(A\ot B)$ where 
$p:Z_{J,D}@>>>\text{point}$ is the obvious map. As in 32.23(b) we have
$$(A:B)=\t(A*(\part_!B).\tag f$$
From definitions we have
$$(A:B)=(\part_!A:\part_!B).\tag g$$
If $A\in\cd(Z_{J,D})$, $B\in\cd(Z_{\e(J),D'}),C\in\cd(Z_{\e'\e(J),D\i D'{}\i})$ then
from (d) we have $\t((A*B)*C)=\t(A*(B*C))$.
Using (f) we rewrite this as $(A*B:\part_!C)=(A:\part_!(B*C))$ or, using (g), as:
$$(A*B:\part_!C)=(\part_!A:B*C).\tag h$$
For $J\sub J',\fc,\fd$ as above and $A\in\cd(Z_{J,D}),B\in\cd(Z_{J',D})$ we show:
$$(B:f_{J,J'}A)=(e_{J,J'}B:A)\in\cd(\text{point}).\tag i$$
Let $p_1,p_2,p_3$ be the obvious maps from $Z_{J',D},Z_{J,D},Z_{J,J',D}$ to the point. We
must show that $p_{1!}(B\ot\fd_!\fc^*A)=p_{2!}(\fc_!\fd^*B\ot A)$. Since 
$p_{1!}\fd_!=p_{3!}=p_{2!}\fc_!$, both sides are equal to $p_{3!}(\fd^*B\ot\fc^*A)$. This
proves (i).

\subhead 36.5\endsubhead
For $J\sub J'\sub\II$ and $w\in{}^{\e(J)}\WW,w'\in{}^{\e(J')}\WW$ we have a diagram 
$${}^wZ_{J,D}@<\fc_{w,w'}<<{}^{w,w'}Z_{J,J',D}@>\fd_{w,w'}>>{}^{w'}Z_{J',D}$$
where 
$${}^{w,w'}Z_{J,J',D}=\{(P,gU_Q)\in Z_{J,J',D};
(P,gU_P)\in{}^wZ_{J,D},(Q,gU_Q)\in{}^{w'}Z_{J',D}\}$$
and $\fc_{w,w'},\fd_{w,w'}$ are the restrictions of $\fc,\fd$ in 36.4.
Consider the functor 

$\ff_{J,w,J',w'}:\cd({}^wZ_{J,D})@>>>\cd({}^{w'}Z_{J',D}),A\m\fd_{w,w'!}\fc_{w,w'}^*A$. 
\nl
From definitions we see that, for $A\in\cd({}^wZ_{J,D})$, we have
$$\ff_{J,w,J',w'}(A)=i_{J',w'}^*\ff_{J,J'}(i_{J,w!}A).\tag a$$

\subhead 36.6\endsubhead
Let $J\sub\II$ and let $w\in{}^{\e(J)}\WW$. We have a sequence of affine space bundles
$${}^wZ_{J,D}@>\vt_{J,w}>>{}^wZ_{J_1,D}@>\vt_{J_1,w}>>{}^wZ_{J_2,D}@>\vt_{J_2,w}>>\do\tag a
$$
where $\vt_{J,w}$ is as in 36.1, $\vt_{J_1,w}$ is the analogous map with $J$ replaced by
$J_1$ (we set $(J_1)_1=J_2$), $\vt_{J_2,w}$ is the analogous map with $J$ replaced by $J_2$
(we set $(J_2)_1=J_3$), etc. We have $J\sps J_1\sps J_2\sps\do$. Let $J_\iy=J_r$ for large
$r$. Since $\vt_{J_r,w}=1$ for large $r$, the composition of the maps (a) is a well defined
map $\un\vt:{}^wZ_{J,D}@>>>{}^wZ_{J_\iy,D}$. For any $(P,gU_P)\in{}^wZ_{J,D}$ we have
$$\vt_{J,w}(P,gU_P)=(P_1,gU_P),\vt_{J_1,w}(P_1,gU_{P_1})=(P_2,gU_{P_2}),\do$$
where $P\sps P_1\sps P_2\sps\do$. Let $P_\iy=P_r$ for large $r$. We show:

(b) {\it The map $\a:{}^{w,w}Z_{J_\iy,J,D}@>>>{}^wZ_{J,D}$, $(R,gU_P)\m(P,gU_P)$, is an 
isomorphism.}
\nl
We show only that $\a$ is a bijection. Let $(P,gU_P)\in{}^wZ_{J,D}$. By 
\cite{\PCS, 4.14(b)} there exists $(B,h)\in\cb\T gU_P$ such that $B\sub P$ and 
$\po(hBh\i,B)=w$. We have also $\po(gBg\i,B)=w$. Let $R=Q_{J_\iy,B}$. By 
\cite{\PCS, 4.14(a)} we have $(R,gU_R)\in{}^wZ_{J_\iy,D}$. We have 
$(R,gU_P)\in{}^{w,w}Z_{J_\iy,J,D}$ and $\a(R,gU_P)=(P,gU_P)$. Thus $\a$ is surjective.

We show that $\a$ is injective. Let $(R,gU_P)\in{}^{w,w}Z_{J_\iy,J,D}$. Define
$P_1,P_2,\do,P_\iy$ in terms of $(P,gU_P)$ as above. It is enough to show that $R=P_\iy$.
Let $P'=gPg\i,R'=gRg\i$. Since $(R,gU_R)\in{}^wZ_{J_\iy,D}$ and $\e\i(\Ad(w)J_\iy)=J_\iy$,
we see that $R,R'$ have a common Levi and $\po(R',R)=w$. Hence if $B'$ is any Borel of 
$R'$, there exists a Borel $B$ of $R$ such that $\po(B',B)=w$. If $w_0=\po(P',P)$, we have
$w=w_0x$ where $x\in\WW_J$ (see \cite{\PCS, 2.1(b)}) and $l(w_0)+l(x)=l(w)$. Hence we can 
find a Borel $B_1$ of $G^0$ such that $\po(B',B_1)=w_0$, $\po(B_1,B)=x$. Since 
$B\sub R\sub P$ we have $B_1\sub P$. Since $\po(B',B_1)=\po(P',P)=w_0$, we have 
$B'\sub(P'\cap P)U_{P'}$ (see \cite{\PCS, 2.7}). Since this holds for any Borel $B'$ of 
$R'$ and $R'$ is the union of its Borels, it follows that $R'\sub(P'\cap P)U_{P'}$. Hence 
$R=g\i R'g\sub(g\i Pg\cap P)U_P=P_1$. Hence $(R,gU_{P_1})\in{}^{w,w}Z_{J_\iy,J_1,D}$ is 
well defined. Repeating the previous argument for $(R,gU_{P_1})$ instead of $(R,gU_P)$ we 
see that $R\sub P_2$. Continuing in this way we obtain $R\sub P_r$ for any $r\ge0$. In 
particular $R\sub P_\iy$. Since $R,P_\iy$ are parabolics of the same type, we see that 
$R=P_\iy$. This proves (b).

Next we show:

(c) {\it If $y\in{}^{\e(J)}\WW$, $y\ne w$, then ${}^{w,y}Z_{J_\iy,J,D}=\em$.}
\nl
Assume that $(R,gU_P)\in{}^{w,y}Z_{J_\iy,J,D}$. We have $(R,gU_R)\in{}^wZ_{J_\iy,D}$ and by
\cite{\PCS, 4.14(b)} there exists $(B,h)\in\cb\T gU_R$ such that $B\sub R$ and 
$\po(hBh\i,B)=w$. We have also $\po(gBg\i,B)=w$. Since $B\sub P$ and 
$w\in{}^{\e(J)}\WW$, we see, using \cite{\PCS, 4.14(a)}, that $(P,gU_P)\in{}^wZ_{J,D}$. 
Since $(P,gU_P)\in{}^yZ_{J,D}$ and ${}^yZ_{J,D},{}^wZ_{J,D}$ are disjoint for $y\ne w$, we
have a contradiction. This proves (c).

Let $A\in{}^w\hZ_{J,D}$. From definitions we have $A=\un\vt^\bst(A')$ with
$A'\in{}^w\hZ_{J_\iy,D}$. We show:

(d) {\it $\ff_{J_\iy,w,J,w}(A')\cong A[-\d]$ for some integer $\d$.}

(e) {\it $\ff_{J_\iy,w,J,y}(A')=0$ for any $y\in{}^{\e(J)}\WW$, $y\ne w$.}
\nl
Now (e) follows immediately from (c). To prove (d) we recall that
$\ff_{J_\iy,w,J,w}(A')=\a_!\ta^*(A')$ where
$\ta:{}^{w,w}Z_{J_\iy,J,D}@>>>{}^wZ_{J_\iy,D}$ is given by $(R,gU_P)\m(R,gU_R)$. Define
$e:{}^wZ_{J,D}@>>>{}^{w,w}Z_{J_\iy,J,D}$ by $e(P,gU_P)=(P_\iy,gU_P)$ with $P_\iy$ as above.
Clearly, $\ta e=\un\vt$, $\a e=1$. Since $\a$ is an isomorphism, we see that
$$\a_!\ta^*(A')=e^*\ta^*(A')=\un\vt^*(A')=\un\vt^\bst(A')[-\d]=A[-\d]$$
where $\d$ is the dimension of any fibre of $\un\vt$. This proves (d).

\subhead 36.7\endsubhead
Let $L\sub J\sub\II$ and let $w\in{}^{\e(J)}\WW,y\in{}^{\e(L)}\WW$. Assume that
$\e\i(\Ad(w)J)=J$, $\e\i(\Ad(y)L)=L$. From definitions we have 
$$\align{}^{y,w}Z_{L,J,D}=&\{(R,gU_P);R\in\cp_L,P=Q_{J,R},gU_P\in D/U_P,\\&
\po(gRg\i,R)=y,\po(gPg\i,P)=w\}.\endalign$$
We show:

(a) ${}^{y,w}Z_{L,J,D}=\em$ unless $y\in\WW_{\e(J)}w=w\WW_J$ that is,
$w=\min(\WW_{\e(J)}y\WW_J)$.
\nl
Assume that $(R,gU_P)\in{}^{y,w}Z_{L,J,D}$. Let $B$ be a Borel of $gRg\i$. Then
$B\sub gPg\i$. In our case $gPg\i,P$ have a common Levi; hence there exists a Borel $B'$ of
$P$ such that $\po(B,B')=\po(gPg\i,P)=w$. Similarly, since $gRg\i,R$ have a common Levi, 
there exists a Borel $B''$ of $R$ such that $\po(B,B'')=\po(gRg\i,R)=y$. Let 
$z=\po(B',B'')$. Since $B'\sub P,B''\sub P$, we have $z\in\WW_J$. Since $w\in\WW^J$ we have
$l(wz)=l(w)+l(z)$ hence from $\po(B,B')=w,\po(B',B'')=z$ we deduce $\po(B,B'')=wz$. Thus 
$y=wz$ and (a) follows.

Now for $J\sub\II$, $w\in{}^{\e(J)}\WW,w'\in{}^{\e'\e(J)}\WW,y\in{}^{\e'\e(J)}\WW$, the
diagram 

$Z_{J,D}\T Z_{\e(J),D'}@<b_1<<Z_0@>b_2>>Z_{J,D'D}$
\nl
(see 32.5) restricts to a diagram 
$${}^wZ_{J,D}\T{}^{w'}Z_{\e(J),D'}@<b_{w,w',1}<<{}^{w,w',y}Z_0@>b_{w,w',2}>>{}^yZ_{J,D'D}$$
where ${}^{w,w',y}Z_0=b_1\i({}^wZ_{J,D}\T {}^{w'}Z_{\e(J),D'})\cap b_2\i({}^yZ_{J,D'D})$.
Define a bifunctor $\cd({}^wZ_{J,D})\T\cd({}^{w'}Z_{\e(J),D'})@>>>\cd({}^yZ_{J,D'D})$ by
$$A,B\m A*_yB=b_{w,w',2!}b_{w,w',1}^*(A\bxt B).$$
Let $K\in\cd({}^wZ_{J,D}),K'\in\cd({}^{w'}Z_{\e(J),D'})$. From definitions we see that
$$i_{J,y}^*((i_{J,w!}K)*(i_{\e(J),w'!}K'))= K*_yK'.\tag b$$

\subhead 36.8\endsubhead
In the remainder of this section we assume that $\kk$ is an algebraic closure of a finite
field. Let $J\sub\II$. As in 31.2 let $\ca=\ZZ[v,v\i]$, $v$ an indeterminate. Let 
$\fK(Z_{J,D})$ (resp. $\fK'(Z_{J,D})$) be the free $\ca$-module with basis given by the 
elements in $\ci_{J,D}$ (resp. in $\sqc_{w\in{}^{\e(J)}\WW}\ci_{J,w,D}$). We have
$\fK'(Z_{J,D})=\op_{w\in{}^{\e(J)}\WW}{}^w\fK'(Z_{J,D})$ where ${}^w\fK'(Z_{J,D})$ is the
$\ca$-submodule of $\fK'(Z_{J,D})$ generated by $\ci_{J,w,D}$. For $w\in{}^{\e(J)}\WW$, 

$i_{J,w}^*:\cd(Z_{J,D})@>>>\cd({}^wZ_{J,D})$,
$i_{J,w!}^*:\cd({}^wZ_{J,D})@>>>\cd(Z_{J,D})$
\nl
restrict to functors

$\cd^{cs}(Z_{J,D})@>>>\cd^{cs}({}^wZ_{J,D})$, 
$\cd^{cs}({}^wZ_{J,D})@>>>\cd^{cs}(Z_{J,D})$,
\nl
see 36.3(a), 36.3(f). These functors, or rather their mixed analogues, are compatible with
the triangular structures hence they induces $\ca$-linear maps 

$\ph_{J,w}:\fK(Z_{J,D})@>>>{}^w\fK'(Z_{J,D})$,
$\ph'_{J,w}:{}^w\fK'(Z_{J,D})@>>>\fK(Z_{J,D})$.
\nl
Clearly, $\ph_{J,w}\ph'_{J,w}=1:{}^w\fK'(Z_{J,D})@>>>{}^w\fK'(Z_{J,D})$. For 
$A\in\ci_{J,D}$ regarded as a pure complex of weight $0$, we have 
$$\ph_{J,w}(A)=\sum_{A_1\in\ci_{J,w,D}}\sum_{j,h}
(-1)^jv^h(\text{multiplicity of $A_1$ in }{}^pH^j(i_{J,w}^*A)_h)A_1$$
(the subscript $h$ means "the subquotient of pure weight $h$"). Define an $\ca$-linear map 
$\ph:\fK(Z_{J,D})@>>>\fK'(Z_{J,D})$ by $\ph=\op_{w\in{}^{\e(J)}\WW}\ph_{J,w}$. 

From 36.3(b),(d) we see that the matrix of $\ph$ is square and upper triangular (with $1$
on diagonal) with respect to a suitable order on $\ci_{J,D}$. In particular,

(a) {\it $\ph:\fK(Z_{J,D})@>>>\fK'(Z_{J,D})$ is an isomorphism.}
\nl
Note that the inverse isomorphism restricted to ${}^w\fK'(Z_{J,D})$ is just $\ph'_{J,w}$.
Clearly,

(b) $\fK(Z_{J,D})=\op_{w\in{}^{\e(J)}\WW}{}^w\fK(Z_{J,D})$,
\nl
where ${}^w\fK(Z_{J,D})=\ph\i({}^w\fK'(Z_{J,D})$ for $w\in{}^{\e(J)}\WW$.

We have a partition $\ci_{J,D}=\sqc_{\fk\in\WW\bsl\ufs(\TT)}\ci_{J,D}^\fk$ where 
$\WW\bsl\ufs(\TT)$ is the set of $\WW$-orbits on $\ufs(\TT)$ and $\ci_{J,D}^\fk$ consists 
of those $A\in\ci_{J,D}$ such that $A\in\hZ_{J,D}^\cl$ for some $\cl\in\fs(\TT)$ whose 
isomorphism class is in $\fk$. We have 

(c) $\fK(Z_{J,D})=\op_{\fk\in\WW\bsl\ufs(\TT)}\fK^\fk(Z_{J,D})$
\nl
where $\fK^\fk(Z_{J,D})$ is the $\ca$-submodule of $\fK(Z_{J,D})$ generated by 
$\ci_{J,D}^\fk$. From definitions we see that each ${}^w\fK'(Z_{J,D})$ and each 
${}^w\fK(Z_{J,D})$ has a natural direct sum decomposition indexed by $\WW\bsl\ufs(\TT)$ 
analogous to (c). It follows that the decompositions (b),(c) are compatible in the sense 
that

(d) $\fK(Z_{J,D})=\op_{w,\fk}({}^w\fK(Z_{J,D})\cap\fK^\fk(Z_{J,D}))$.

\subhead 36.9\endsubhead
Let $J\sub J'\sub\II$. As in \cite{\PCS, 6.4, 6.7(b)} we see that $\ff_{J,J'},\fe_{J,J'}$
restrict to functors $\ff_{J,J'}:\cd^{cs}(Z_{J,D})@>>>\cd^{cs}(Z_{J',D})$,
$\fe_{J,J'}:\cd^{cs}(Z_{J',D})@>>>\cd^{cs}(Z_{J,D})$, denoted again by 
$\ff_{J,J'},\fe_{J,J'}$. These functors or, rather, their mixed analogues, are compatible 
with the triangulated structures hence induce $\ca$-linear maps 
$\fK(Z_{J,D})@>>>\fK(Z_{J',D})$, $\fK(Z_{J',D})@>>>\fK(Z_{J,D})$
denoted again by $\ff_{J,J'},\fe_{J,J'}$; the identities 36.4(b) continue to hold for these
linear maps. Moreover, Verdier duality induces a group homomorphism 
$\fK(Z_{J,D})@>>>\fK(Z_{J,D})$ (denoted again by $\fD$) and from 36.4(a) we deduce
$$\fD\ff_{J,J'}=v^{-2a}\ff_{J,J'}\fD:\fK(Z_{J,D})@>>>\fK(Z_{J',D})\tag a$$ 
with $a$ as in 36.4(a). Now for any $A\in\ci_{J,D}$ we have
$$\fD(v^mA)=v^{-m}A^*\tag b$$
where $A^*\in\ci_{J,D}$ is isomorphic to the Verdier dual of $A$. If $J\sub J'$, 
$\fc,\fd,a$ are as in 36.4, then for $A\in\ci_{J,D},B\in\ci_{J',D}$ (regarded as pure 
perverse sheaves of weight $0$) we have 
$$\ff_{J,J'}A=\sum_{A_1\in\ci_{J',D}}\sum_{j,h}(-1)^jv^h(\text{multiplicity of $A_1$ in }
{}^pH^j(\fd_!\fc^*A)_h)A_1,$$
$$\fe_{J,J'}B=\sum_{B_1\in\ci_{J,D}}\sum_{j,h}(-1)^jv^h(\text{multiplicity of $B_1$ in }
{}^pH^j(\fc_!\fd^*B)_h)B_1$$
(the subscript $h$ means "the subquotient of pure weight $h$"). Since $\fc$ is an affine
space bundle with fibres of dimension $a$, $\fc^*A[a]$ is a pure perverse sheaf of weight
$a$. By \cite{\DE, 6.2.6} applied to the proper morphism $\fd$, $\fd_!\fc^*A[a]$ is a pure
complex of weight $a$. From \cite{\BBD, 5.4.1}, we see that
${}^pH^j(\fd_!\fc^*A)={}^pH^{j-a}(\fd_!\fc^*A[a])$ is pure of weight $(j-a)+a=j$. Hence,
writing $\ff_{J,J'}A=\sum_{A_1\in\ci_{J',D}}x_{A,A_1}A_1, x_{A,A_1}\in\ca$ we have
$$x_{A,A_1}=\sum_j(-v)^j(\text{multiplicity of $A_1$ in }{}^pH^{j-a}(\fd_!\fc^*A[a])).$$
Let $\x\m\bar\x$ be the group homomorphism $\ca@>>>\ca$ given by $v^m\m v^{-m}$ for 
$m\in\ZZ$. Since $\fc^*A[a]$ is a pure perverse sheaf we have (using the relative hard 
Lefschetz theorem \cite{\BBD, 5.4.10}):
$$\align&\bar x_{A,A_1}
=\sum_j(-v)^{-j}(\text{multiplicity of $A_1$ in }{}^pH^{j-a}(\fd_!\fc^*A[a]))\\&
=\sum_j(-v)^{-j}(\text{multiplicity of $A_1$ in }{}^pH^{-j+a}(\fd_!\fc^*A[a]))\\&
=\sum_{j'}(-v)^{j'-2a}(\text{multiplicity of $A_1$ in }{}^pH^{j'-a}(\fd_!\fc^*A[a]))
=v^{-2a}x_{A,A_1}.\endalign$$
Define a group homomorphism $\b_J:\fK(Z_{J,D})@>>>\fK(Z_{J,D})$ by $\b_J(v^mA)=v^{-m}A$ for
$A\in\ci_{J,D},m\in\ZZ$. We see that $\b_{J'}(\ff_{J,J'}A)=v^{-2a}\ff_{J,J'}A$ for any 
$A\in\ci_{J,D},m\in\ZZ$. Hence for any $\x\in\fK(Z_{J,D})$ we have
$$\b_{J'}(\ff_{J,J'}\x)=v^{-2a}\ff_{J,J'}(\b_J(\x)).\tag c$$
Define $\ti{\fD}:\fK(Z_{J,D})@>>>\fK(Z_{J,D})$ by $\ti{\fD}A=\fD\b_J(A)=\b_J\fD(A)$. 
Clearly, $\ti{\fD}$ is $\ca$-linear. Note that 

(d) {\it the maps $\ff_{J,J'},\fe_{J,J'},\b_J$ are compatible with the decompositions of
type 36.8(c) while $\fD,\ti{\fD}$ map the summand corresponding to $\fk$ in 36.8(c) onto 
the summand corresponding to $\che\fk$. (Here $\che\fk$ is the image of $\fk$ under 
$\cl\m\che\cl$, see 28.18.)}
\nl
From 32.21 we see that convolution (see 32.5) restricts to a bifunctor

$\cd^{cs}(Z_{J,D})\T\cd^{cs}(Z_{\e(J),D'})@>>>\cd^{cs}(Z_{J,D'D})$.
\nl
This induces an $\ca$-bilinear pairing 
$\fK(Z_{J,D})\T\fK(Z_{\e(J),D'})@>>>\fK(Z_{J,D'D})$. It sends 
$A\in\ci_{J,D},B\in\ci_{\e(J),D'}$ (regarded as pure perverse sheaves of weight $0$) to
$$\sum_{A_1\in\ci_{J,D'D}}\sum_{j,h}(-1)^jv^h
(\text{multiplicity of $A_1$ in }{}^pH^j(A*B)_h)A_1.$$ 
This bilinear pairing is denoted again by $A,B\m A*B$. If $\fk,\fk'\in\WW\bsl\ufs(\TT)$ and
$\x\in\ufK^\fk(Z_{J,D})$, $\x'{}^{\fk'}\in\fK(Z_{\e(J),D'})$, then 
$\x*\x'\in\fK^\fk(Z_{J,D'D})$ if $\fk'=\uD(\fk)$ and $\x*\x'=0$ if $\fk'\ne\uD(\fk)$ (see 
32.6(a),(b)).

From 28.19 we see that $\part_!:\cd(Z_{J,D})@>>>\cd(Z_{\e(J),D\i})$ restricts to an
equivalence of categories $\cd^{cs}(Z_{J,D})@>>>\cd^{cs}(Z_{\e(J),D\i})$; this induces an 
$\ca$-linear isomorphism $\fK(Z_{J,D})@>\si>>\fK(Z_{\e(J),D\i})$ denoted again by $\part$. 
Now the composition $\fK(Z_{J,D})@>\part>>\fK(Z_{\e(J),D\i})@>\part>>\fK(Z_{J,D})$ is the 
identity map. 

In the setup of 36.5,  $\ff_{J,w,J',w'}$ restricts to a functor
$\ff_{J,w,J',w'}:\cd^{cs}({}^wZ_{J,D})@>>>\cd^{cs}({}^{w'}Z_{J',D})$. (We use the analogous
statement for $\ff_{J,J'}$, as above, together with 36.5(a), 36.3(a), 36.3(f).) Hence 
there is an induced $\ca$-linear map ${}^w\fK'(Z_{J,D})@>>>{}^{w'}\fK'(Z_{J',D})$ denoted 
again by $\ff_{J,w,J',w'}$. Let $\ff'_{J,J'}:\fK'(Z_{J,D})@>>>\fK'(Z_{J',D})$ be the unique
$\ca$-linear map such that $pr_{w'}\ff'_{J,J'}pr_w=\ff_{J,w,J',w'}$ for any $w,w'$ as 
above. (Here $pr_w$ is the projection of $\fK'(Z_{J',D})$ onto the direct summand 
${}^w\fK'(Z_{J,D})$.) From 36.5(a) we see that for $x\in{}^w\fK'(Z_{J',D})$ we have
$$\ff'_{J,J'}=\ph\ff_{J,J'}\ph\i.\tag e$$
We also have
$$\part\ff_{J,J'}=\ff_{\e(J),\e(J')}\part,\part\fe_{J,J'}=\fe_{\e(J),\e(J')}\part.\tag f
$$
Now the functor $\t:\cd^{cs}(Z_{J,G^0})@>>>\cd(\text{point})$ (see 36.4) or rather its 
mixed analogue induces an $\ca$-linear map $\fK(Z_{J,G^0})@>>>\ca$; it sends 
$C\in\ci_{J,G^0}$ regarded as a pure complex of weight $0$ to 

$\sum_{j,h}(-1)^jv^h\dim(H^j_c(\cp_J,i^*C)_h)\in\ca$
\nl
where $i:\cp_J@>>>Z_{J,G^0}$ is $P\m(P,U_P)$ and $H^j_c(\cp_J,i^*C)_h$ denotes the part of
weight $h$ of the mixed vector space $H^j_c(\cp_J,i^*C)$. This linear map is denoted again
by $\t$. 

Similarly, the bifunctor $(:):\cd^{cs}(Z_{J,D})\T\cd^{cs}(Z_{J,D})@>>>\cd(\text{point})$ 
(see 36.4) or rather its mixed analogue induces an $\ca$-bilinear pairing 
$\fK(Z_{J,D})\T\fK(Z_{J,D})@>>>\ca$; it sends $A,B\in\ci_{J,D}$ regarded as pure complexes 
of weight $0$ to 

$\sum_{j,h}(-1)^jv^h\dim(H^j_c(Z_{J,D},A\ot B)_h)\in\ca$.
\nl
(Here $H^j_c(Z_{J,D},A\ot B)_h$ denotes the part of weight $h$ of the mixed vector space
$H^j_c(Z_{J,D},A\ot B)$.) This bilinear pairing is denoted again by $(:)$. From 
\cite{\BBD, 5.1.14} we see that $(A:B)\in\ZZ[v\i]$ and from \cite{\CS, II,7.4} and its 
proof we see that the constant term of $(A:B)$ is $\d_{A,B^*}$. Thus,
$$(A:B)\in\d_{A,B^*}+v\i\ZZ[v\i].\tag g$$
We show:

(h) {\it if $\x\in\fK(Z_{J,D}),(\x:\ti{\fD}(\x)))=0$ then $\x=0$.}
\nl
Assume that $\x\ne0$. We have $\x=\sum_{A\in\ci_{J,D},m\in\ZZ}g_{A,m}v^mA$ where 
$g_{A,m}\in\ZZ$ is zero for all but finitely many $A,m$. We can find $e\in\ZZ$ such that
$g_{A,e}\ne0$ for some $A\in\ci_{J,D}$ and $g_{A,m}=0$ for all $m>e$ and all 
$A\in\ci_{J,D}$. Then

$(\x:\ti{\fD}(\x)))=\sum_{A_1,A,m_1,m}g_{A_1,m_1}g_{A,m}v^{m_1+m}(A_1:A^*)$.
\nl
By (g) this equals $\sum_Ag_{A,e}^2v^{2e}+\text{an element in }v^{2e-1}\ZZ[v\i]$. This is 
non-zero since $\sum_Ag_{A,e}^2\in\ZZ_{>0}$. This proves (h).

We show:

(i) {\it If $\fk,\fk'\in\WW\bsl\ufs(\TT)$ and $\fk'\ne\che\fk$ then 
$(\fK^\fk(Z_{J,D}):\fK^{\fk'}(Z_{J,D}))=0$.}
\nl
Let $A\in\hZ_{J,D}^\cl,B\in\hZ_{J,D}^{\cl_1}$ where $\cl,\cl_1\in\fs(\TT)$. It is enough
to show that:

{\it if $H^j_c(Z_{J,D},A\ot B)\ne0$ for some $j$ then the isomorphism class of $\cl_1$ is 
in the $\WW$-orbit of the isomorphism class of $\che\cl$.}
\nl
From our assumption we have $A*\part_!(B)\ne0$. We have $A\dsv\bK^{\bs,\cl}_{J,D}$ as in 
28.13(v) and similarly 
$B\dsv\bK^{\bs',\cl_1}_{J,D}$. Hence $\part_!B\dsv\bK^{\bs',\cl''}_{\e(J),D\i}$ where
$\cl''=(\uD\i)^*\che\cl$, see 28.19. We then have
$\bK^{\bs,\cl}_{J,D}*\bK^{\bs',\cl''}_{\e(J),D\i}\ne0$. Using 32.6(a) we see that $\cl$,
$\cl_1$ have the required property. This proves (i).

We show:

(j) {\it If $w,w'\in\WW$ and $w\ne w'$ then $({}^w\fK(Z_{J,D}):{}^{w'}\fK(Z_{J,D}))=0$.}
\nl
It is enough to show: if $A\in\cd^{cs}({}^wZ_{J,D}),A'\in\cd^{cs}({}^{w'}Z_{J,D})$ then
$p_!(i_{J,w!}A\ot i_{J,w'!}A')=0$ with $p$ as in 36.4. It is enough to show that
$(i_{J,w}\T i_{J,w'})_!f^*(A\bxt A')=0$ where $f:Z_{J,D}@>>>Z_{J,D}\T Z_{J,D}$ is the 
diagonal map. This follows from the fact that ${}^wZ_{J,D}\cap({}^{w'}Z_{J,D})=\em$.

We show:

(k) {\it If $w\in{}^{\e(J)}\WW$, $w\ne1$ then $\t({}^w\fK(Z_{J,G^0}))=0$.}
\nl
Since ${}^wZ_{J,G^0},{}^1Z_{J,G^0}$ are disjoint, we have $i_{J,1}^*i_{J,w!}=0$ hence 
$i^*i_{J,w!}=0$ with $i$ as in 36.4; (k) follows.

From the definitions we have

(l) $\part({}^w\fK(Z_{J,D}))={}^{w\i}\fK(Z_{\e(J),D\i})$
\nl
for any $w\in{}^{\e(J)}\WW$ such that $\e\i(\Ad(w)J)=J$. (For such $w$ we have
$w\i\in{}^J\WW$.)

\subhead 36.10\endsubhead
Let 

$\tfK(Z_{J,D})=\sum_{L;L\subsetneq J}\ff_{L,J}(\fK(Z_{L,D}))$,
$\fK_0=\op_{w\in{}^{\e(J)}\WW;\e\i(\Ad(w)J)\ne J}{}^w\fK(Z_{J,D})$.
\nl
We show:
$$\tfK(Z_{J,D})=\fK_0\op\op\Sb w\in{}^{\e(J)}\WW\\ \e\i(\Ad(w)J)=J\endSb
(\tfK(Z_{J,D})\cap{}^w\fK(Z_{J,D})).\tag a$$
Setting 

$\tfK'(Z_{J,D})=\sum_{L;L\subsetneq J}\ff'_{L,J}(\fK'(Z_{L,D}))$,
$\fK'_0=\op_{w\in{}^{\e(J)}\WW;\e\i(\Ad(w)J)\ne J}{}^w\fK'(Z_{J,D})$ 
\nl
and using 36.9(e), we see that it is enough to show:
$$\tfK'(Z_{J,D})=\fK_0\op\op\Sb w\in{}^{\e(J)}\WW\\ \e\i(\Ad(w)J)=J\eSb
(\tfK'(Z_{J,D})\cap{}^w\fK'(Z_{J,D})).\tag b$$
From 36.6(d), 36.6(e) we see that, if $w\in{}^{\e(J)}\WW,\e\i(\Ad(w)J)\ne J$, then
$A=\ff'_{J_\iy,J}A'$ for some $A'\in{}^y\fK'(Z_{L',D})$ and $J_\iy\ne J$. Hence

(c) $\fK'_0\sub\tfK'(Z_{J,D})$.
\nl
Using again 36.6(d), 36.6(e) and also 36.4(b) (or rather its analogue for $\ff'_{J,J'}$)
and (c), we see that
$$\align&\tfK'(Z_{J,D})=\sum\Sb L;L\subsetneq J\\y\in{}^{\e(L)}\WW\\ \e\i(\Ad(y)L)=L\eSb
\ff'_{L,J}({}^y\fK'(Z_{L,D}))\\&=\fK'_0+\sum\Sb L;L\subsetneq J\\y\in{}^{\e(L)}\WW\\
\e\i(\Ad(y)L)=L\eSb\ff'_{L,J}({}^y\fK'(Z_{L,D})).\endalign$$
Thus, (b) holds and (a) holds.

Let 
$$\bfK(Z_{J,D})=\fK(Z_{J,D})/\tfK(Z_{J,D}).$$
For any $w\in{}^{\e(J)}\WW$ such that $\e\i(\Ad(w)J)=J$, let ${}^w\bfK(Z_{J,D})$ be the 
image of ${}^w\fK(Z_{J,D})$ under the obvious map $\fK(Z_{J,D})@>>>\bfK(Z_{J,D})$. From 
(a) we see that
$$\bfK(Z_{J,D})=\op\Sb w\in{}^{\e(J)}\WW\\ \e\i(\Ad(w)J)=J\eSb{}^w\bfK(Z_{J,D}).\tag d$$
For any $\fk\in\WW\bsl\ufs(\TT)$ let $\bfK^\fk(Z_{J,D})$ be the image of 
$\fK^\fk(Z_{J,D})$ under the obvious map $\fK(Z_{J,D})@>>>\bfK(Z_{J,D})$. From 36.9(d) we
see that $\tfK(Z_{J,D})$ is compatible with the decomposition 36.8(c). It follows that 
$$\bfK(Z_{J,D})=\op_{\fk\in\WW\bsl\ufs(\TT)}\bfK^\fk(Z_{J,D}).\tag e$$
Using 36.8(d) we see that that the decompositions (d),(e) are compatible in the sense that
$$\bfK(Z_{J,D})=\op_{w,\fk}({}^w\bfK(Z_{J,D})\cap\bfK^\fk(Z_{J,D})).\tag f$$
Using 36.9(a), 36.9(c), we see that

(g) {\it $\fD,\b_J,\ti{\fD}:\fK(Z_{J,D})@>>>\fK(Z_{J,D})$ map $\tfK(Z_{J,D})$ into itself 
hence induces group homomorphisms $\bfK(Z_{J,D})@>>>\bfK(Z_{J,D})$ denoted again by 
$\fD,\b_J,\ti{\fD}$.}

\subhead 36.11\endsubhead
For $J'\sub J\sub\II$ and $x\in\fK(Z_{J',D}),x'\in\fK(Z_{\e(J),D'})$, we show:
$$(\ff_{J',J}x)*x'=\ff_{J',J}(x*\fe_{\e(J'),\e(J)}x').\tag a$$
Using successively 36.4(h), 36.9(f), 36.4(i), 36.4(c), 36.9(f), 36.4(h), 36.4(i), we see
that, for any $x''\in\fK(Z_{J,D'D})$, we have
$$\align&((\ff_{J',J}x)*x':x'')=((\part\ff_{J',J}x):x'*\part x'')\\&=
(\ff_{\e(J'),\e(J)}\part(x):x'*(\part x''))=(\part x:\fe_{\e(J'),\e(J)}(x'*(\part x'')))
\\&=(\part x:\fe_{\e(J'),\e(J)}(x')*\fe_{\e(J'),\e(J)}(\part x''))) 
=(\part x:\fe_{\e(J'),\e(J)}x'*\part\fe_{J',J}x'') \\&
=(x*\fe_{\e(J'),\e(J)}x':\fe_{J',J}x'')=(\ff_{J',J}(x*\fe_{\e(J'),\e(J)}x'):x'').\endalign
$$
Thus, if $\x=(\ff_{J',J}x)*x'-\ff_{J',J}(x*\fe_{\e(J'),\e(J)}x')$ then $(\x:x'')=0$ for any
$x''\in\fK(Z_{J,D'D})$. In particular, $(\x:\ti{\fD}(\x))=0$. Using 36.9(h) we see that
$\x=0$. This proves (a). 

A similar argument shows that, if $x\in\fK(Z_{J,D}),x'\in\fK(Z_{\e(J'),D'})$, then:
$$x*(\ff_{\e(J'),\e(J)}x')=\ff_{J',J}(\fe_{J',J}x*x').\tag b$$
From (a),(b) we see that 

$\tfK(Z_{J,D})*\fK(Z_{\e(J),D'})\sub\tfK(Z_{J,D'D})$,
$\fK(Z_{J,D})*\tfK(Z_{\e(J),D'})\sub\tfK(Z_{J,D'D})$.
\nl
It follows that $\fK(Z_{J,D})\T\fK(Z_{\e(J),D'})@>>>\fK(Z_{J,D'D})$, $A,B\m A*B$, induces
an $\ca$-bilinear pairing $\bfK(Z_{J,D})\T\bfK(Z_{\e(J),D'})@>>>\bfK(Z_{J,D'D})$. We 
denote it again by $A,B\m A*B$. The following result relates this bilinear pairing to the
decompositions of type 36.10(d).
\nl
For $w\in{}^{\e(J)}\WW,w'\in{}^{\e'\e(J)}\WW$ such that $\e\i(\Ad(w)J)=J$, 
$\e'{}\i(\Ad(w')\e(J))=\e(J)$, let $X_{w,w'}$ be the set of all $y\in{}^{\e'\e(J)}\WW$ such
that $\e\i\e'{}\i(\Ad(y)J)=J$ and such that for some 
$Q\in\cp_J,Q'\in\cp_{\e(J)},Q''\in\cp_{\e'\e(J)}$ we have

$\po(Q',Q)=w,\po(Q'',Q')=w',\po(Q'',Q)=y$.
\nl
Then
$${}^w\bfK(Z_{J,D})*{}^{w'}\bfK(Z_{\e(J),D'})\sub\op_{y\in X_{w,w'}}{}^y\bfK(Z_{J,D'D}).
\tag c$$
Indeed, using 36.7(a) we see that it is enough to show:

{\it If ${}^{w,w',y}Z_0\ne\em$ (notation of 36.7) then $y\in X_{w,w'}$.}
\nl
This is immediate from definitions.

\subhead 36.12\endsubhead
As a special case of the pairing $A,B\m A*B$ in 36.11 we have an $\ca$-bilinear pairing
$\bfK(Z_{J,G^0})\T\bfK(Z_{J,G^0})@>>>\bfK(Z_{J,G^0})$. This defines an associative algebra
structure on $\bfK(Z_{J,G^0})$ (not necessarily with $1$). 

Also as special cases of the pairing $A,B\m A*B$ in 36.11 we have $\ca$-bilinear pairings
$$\bfK(Z_{J,G^0})\T\bfK(Z_{J,D})@>>>\bfK(Z_{J,D}),
\bfK(Z_{J,D})\T\bfK(Z_{\e(J),G^0})@>>>\bfK(Z_{J,D}),$$
which make $\bfK(Z_{J,D})$ into a (not necessarily unital) 
$(\bfK(Z_{J,G^0}),\bfK(Z_{\e(J),G^0}))$-bimodule.

By 36.10(f) we have
$$\bfK^\he(Z_{J,G^0})=\op\Sb w\in{}^{\e(J)}\WW\\ \e\i(\Ad(w)J)=J\eSb
({}^w\bfK(Z_{J,G^0})\cap\bfK^\he(Z_{J,G^0})).$$
From 36.10(g) we see that $\fD,\b_J,\ti{\fD}$ may be regarded as group homomorphisms
$\bfK^\he(Z_{J,G^0})@>>>\bfK^\he(Z_{J,G^0})$.

\subhead 36.13\endsubhead
For any $\ca$-module $V$ we set $\un{V}=\QQ(v)\ot_\ca V$. If $V,V'$ are $\ca$-modules and 
$f:V@>>>V'$ is $\ca$-linear we denote again by $f$ the $\QQ(v)$-linear map 
$\un V@>>>\un V'$ such that $1\ot x\m 1\ot f(x)$ for $x\in V$.

In particular, the $\QQ(v)$-vector spaces $\ufK(Z_{J,D})$, $\un{\tfK}(Z_{J,D})$ and
$\un{\bfK}(Z_{J,D})$ are defined. Note that 

$\un{\tfK}(Z_{J,D})=\sum_{L;L\subsetneq J}\ff_{L,J}(\ufK(Z_{L,D}))$,
$\ufK(Z_{J,D})=\op_{\fk\in\WW\bsl\ufs(\TT)}\ufK^\fk(Z_{J,D})$,
\nl
see 36.8(c). 

The symmetric bilinear form $(:):\fK(Z_{J,D})\T\fK(Z_{J,D})@>>>\ca$ extends uniquely to a 
symmetric $\QQ(v)$-bilinear form $\ufK(Z_{J,D})\T\ufK(Z_{J,D})@>>>\QQ(v)$ denoted again by
$(:)$.

A vector subspace $E$ of $\ufK(Z_{J,D})$ is said to be {\it homogeneous} if 
$E=\sum_{\fk\in\WW\bsl\ufs(\TT)}E^\fk$ where $E^\fk=E\cap\ufK^\fk(Z_{J,D})$. For such $E$ 
we set $E^\pe=\{x\in\ufK(Z_{J,D});(x:E)=0\}$. Then $E^\pe$ is homogeneous. We show that if,
in addition, $E$ is stable under $\ti{\fD}$, then:

(a) $\ufK(Z_{J,D})=E\op E^\pe$.
\nl
We first show that $E\cap E^\pe=0$. Assume that $x\in E\cap E^\pe$. Then $\ti{fD}(x)\in E$
and $(x:\ti{\fD}(x))=0$. We can find $\l\in\ca-\{0\}$ such that $\l x\in\fK(Z_{J,D})$. Then
$(\l x,\ti{\fD}(\l x))=0$. By 36.9(h) we have $\l x=0$. Hence $x=0$ so that 
$E\cap E^\pe=0$. It remains to show that $\dim(E^\pe)^\fk+\dim E^\fk=\dim\ufK^\fk(Z_{J,D})$
for any $\fk$. This is a consequence of the following statement:

(b) {\it for any $\fk$, the form $(:)$ on the finite dimensional vector space
$\ufK^\fk(Z_{J,D})+\ufK^{\che\fk}(Z_{J,D})$ is non-singular.}
\nl
More generally we show that for any finite dimensional subspace $E'$ of $\ufK(Z_{J,D})$ 
which is stable under $\ti{\fD}$, the form $(:)$ on $E'$ is non-singular. Let $x\in E'$ be
such that $(x:x')=0$ for any $x'\in E'$. In particular we have $(x,\ti{\fD}(x))=0$. As 
above, this implies that $x=0$. This proves (b) hence also (a).

We set 
$$\ufK(Z_{J,D})^J=\{\x\in\ufK(Z_{J,D});\fe_{H,J}\x=0\qua\frl H\subsetneq J\}.$$
We show:
$$\un{\tfK}(Z_{J,D})^\pe=\ufK(Z_{J,D})^J.\tag c$$
We have
$$\align&\un{\tfK}(Z_{J,D})=\{x\in\ufK(Z_{J,D});(x:\ff_{H,J}\ufK(Z_{H,D}))=0\qua\frl 
H\subsetneq J\}\\&=\{\x\in\ufK(Z_{J,D});(\fe_{H,J}x:\ufK(Z_{H,D}))=0\qua\frl H\subsetneq J
\}\\&=\{\x\in\ufK(Z_{J,D});\fe_{H,J}x=0\qua\frl H\subsetneq J\},\endalign$$
as required. (We have used that $\ufK(Z_{H,D})^\pe=0$ which follows from (a), applied to 
$H$ instead of $J$.) This proves (c). 

We show:
$$\ufK(Z_{J,D})=\un{\tfK}(Z_{J,D})\op\ufK(Z_{J,D})^J.\tag d$$
This follows from (a),(c) using the fact that $\un{\tfK}(Z_{J,D})$ is stable under 
$\ti{\fD}$, see 36.10(g).

From (d) we see that the second projection $\ufK(Z_{J,D})@>>>\ufK(Z_{J,D})^J$ induces an 
isomorphism $\ufK(Z_{J,D})/\un{\tfK}(Z_{J,D})@>\si>>\ufK(Z_{J,D})^J$ that is, an 
isomorphism 
$$\un{\bfK}(Z_{J,D})@>\si>>\ufK(Z_{J,D})^J.\tag e$$

\subhead 36.14\endsubhead
The $\ca$-algebra structure on $\fK(Z_{J,G^0})$ (resp. $\bfK(Z_{J,G^0})$) given by $*$ 
extends to a $\QQ(v)$-algebra structure on $\ufK(Z_{J,G^0})$ (resp. $\un{\bfK}(Z_{J,G^0})$)
denoted again by $*$. Moreover, $\un{\tfK}(Z_{J,G^0})$ is a two-sided ideal of
$\ufK(Z_{J,G^0})$ (see 36.11). Note that
$\ufK(Z_{J,D})^J$ is also a two-sided ideal of $\ufK(Z_{J,G^0})$. This follows from the 
fact that $\fe_{H,J}:\ufK(Z_{J,G^0})@>>>\ufK(Z_{H,G^0})$ is an algebra homomorphism for any
$H\sub J$ (see 36.4(c)). Since the two summands in the right hand side of 36.13(d) are 
two-sided ideals, they annihilate each other under the product $*$. We see also that the 
isomorphism 36.13(e) respects the algebra structures.

From 36.9(f) we see that the two summands in the right hand side of 36.13(d) are stable 
under $\part:\ufK(Z_{J,G^0})@>>>\ufK(Z_{J,G^0})$ and from 36.4(h) we have
$$(x*y:\part(z))=(\part(x):y*z)\tag a$$
for $x,y,z\in\ufK(Z_{J,G^0})^J$. 

Let $\cw:=\{w\in{}^J\WW;\Ad(w)J=J\}$. We show:
$$\ufK(Z_{J,G^0})^J=\op_{w\in\cw}{}^w\ufK(Z_{J,G^0})^J\tag b$$
where ${}^w\ufK(Z_{J,G^0})^J={}^w\ufK(Z_{J,G^0})\cap\ufK(Z_{J,G^0})^J$. Let 
$x\in\ufK(Z_{J,G^0})^J$. By 36.8(b) we can write uniquely $x=\sum_{w\in{}^J\WW}x_w$ where 
$x_w\in{}^w\ufK(Z_{J,G^0})$. It is enough to show that $x_w\in\ufK(Z_{J,G^0})^J$ (that is,
$(y:x_w)=0$ for any $y\in\un{\tfK}(Z_{J,G^0})$) for all $w$ and $x_w=0$ (that is, 
$(y':x_w)=0$ for any $y'\in\ufK(Z_{J,G^0})$) if $\Ad(w)J\ne J$. For $x_w,y$ as above we 
have $y=\sum_{w'\in{}^J\WW}y_{w'}$ with 
$y_{w'}\in{}^{w'}\ufK(Z_{J,G^0}\cap\un{\tfK}(Z_{J,G^0})$, see 36.10(a). Using 36.9(j) 
twice, we have $(y:x_w)=(y_w:x_w)=(y_w:x)=0$, as required. Now assume that $w\in{}^J\WW$,
$\Ad(w)J\ne J$ and $y'\in\ufK(Z_{J,G^0})$. By 36.8(b) we have 
$y'=\sum_{w'\in{}^J\WW}y'_{w'}$ with $y'_{w'}\in{}^{w'}\ufK(Z_{J,G^0}$; moreover, by
36.10(a) we have $y'_w\in\un{\tfK}(Z_{J,G^0})$. Using 36.9(j) twice, we have 
$(y':x_w)=(y'_w:x_w)=(y'_w:x)=0$, as required. This proves (b).

From definitions we see that 

(c) {\it the decomposition (b) corresponds under 36.13(e) to the decomposition\lb
$\un{\bfK}(Z_{J,G^0})=\op_{w\in\cw}{}^w\un{\bfK}(Z_{J,G^0})$ (see 36.10(d)).}
\nl
Under the isomorphism 36.13(e), the involution 
$\ti{\fD}:\un{\bfK}(Z_{J,D})@>>>\un{\bfK}(Z_{J,D})$ corresponds to an involution
$\ti{\fD}':\ufK(Z_{J,D})^J@>>>\ufK(Z_{J,D})^J$; this is related to the involution 
$\ti{\fD}:\ufK(Z_{J,D})@>>>\ufK(Z_{J,D})$ by
$\ti{\fD}'x=\ti{\fD}x\mod\un{\tfK}(Z_{J,D}$ for $x\in\ufK(Z_{J,D})^J$. Hence for 
$x,x'\in\ufK(Z_{J,D})^J$ we have

$(\ti{\fD}'x:x')=(\ti{\fD}x:x')$.
\nl
Using this and 36.9(h) we see that:

(d) {\it if $x\in\ufK(Z_{J,D})^J$ satisfies $(x:\ti{\fD}x)=0$ then $x=0$.}
\nl
Writing $\he$ instead of $\{\bbq\}\sub\ufs(\TT)$ we see that $\ufK^\he(Z_{J,G^0})$ is a 
subalgebra of $\ufK(Z_{J,G^0})$. Hence 
$$\ufK^J:=\ufK^\he(Z_{J,D})\cap\ufK(Z_{J,D})^J$$
is a subalgebra of $\ufK(Z_{J,D})^J$. From (b) and 36.8(d) we see that
$$\ufK^J=\op_{w\in\cw}{}^w\ufK^J\tag e$$
where ${}^w\ufK^J={}^w\ufK^\he(Z_{J,G^0})\cap\ufK^J$. Now 
$\part:\ufK(Z_{J,G^0})^J@>>>\ufK(Z_{J,G^0})^J$ leaves $\ufK^J$ stable.
Moreover, for $w\in\cw$ we have
$$\part({}^w\ufK^J)={}^{w\i}\ufK^J.\tag f$$
(see 36.9(l)). 

Note that $\cw$ is the same as the 
set of all $w\in\WW$ such that the corresponding permutation of the set of roots leaves 
stable the set of simple roots corresponding to elements of $J$. Hence $\cw$ is a subgroup
of $\WW$.

\subhead 36.15\endsubhead
By 36.4(f) we have $(x:x')=(x':x)=\t(x*\part(x'))$ for $x,x'\in\ufK^J$.
By 36.13(b), the bilinear form $(:)$ is non-degenerate on (the finite dimensional vector 
space) $\ufK^J$. Hence there is a unique vector $x_0\in\ufK^J$ such that
$(x_0:x)=\t(x)$ for all $x\in\ufK^J$. Hence for $x,x'\in\ufK^J$ we have
$(x_0:x*\part(x'))=\t(x*\part(x'))=(x:x')$. Using 36.14(a) we rewrite this as
$(\part(x_0)*x:x')=(x:x')$. (We use also that $\part^2=1$ on $\ufK^J$.) Using the
non-degeneracy of $(:)$ on $\ufK^J$ we deduce $\part(x_0)*x=x$ for all $x\in\ufK^J$. For 
$x,x'\in\ufK^J$ we have also $(\part(x)*x':x_0)=\t(\part(x)*x')=(\part(x):\part(x'))
=(x:x')$ (we use 36.4(g)). Using 36.14(a) we rewrite this as
$(x:x'*\part(x_0))=(x:x')$. Using the
non-degeneracy of $(:)$ on $\ufK^J$ we deduce $x'*\part(x_0)=x'$ for all $x'\in\ufK^J$.
We see that {\it the algebra $\ufK^J$ has a unit element, namely $\bold1=\part(x_0)$.}
By 36.14(e) we have $x_0=\sum_{w\in\cw}x^w_0$ where $x^w_0\in{}^w\ufK^J$. For any 
$x\in\ufK^J$ we have similarly $x=\sum_{w\in\cw}x^w$ where $x^w\in{}^w\ufK^J$. From 
36.9(k) we see that $\t(x)=\t(x^1)$. Hence $(x_0:x^1)=\t(x^1)=\t(x)=(x_0:x)$. Using 
36.9(j) we see that $(x_0:x^1)=(x_0^1:x^1)=(x_0^1:x)$ hence $(x_0^1:x)=(x_0:x)$. Using 
the non-degeneracy of $(:)$ on $\ufK^J$ we deduce $x_0=x_0^1$ that is $x_0\in{}^1\ufK^J$.
Using 36.14(f) we deduce that $\bold1\in{}^1\ufK^J$.

\subhead 36.16\endsubhead
We preserve the setup of 36.14, 36.15. In this subsection we assume that $G=G^0$ is a symplectic 
group $Sp_{2n}(\kk)$ ($n\ge1$) and we describe in this case the structure of the algebra
$\ufK^J$. (The proofs, which depend on results in this and future Sections, 
will be given elsewhere.) We have $\II=\{s_1,s_2,\do,s_n\}$ where $s_is_{i+1}$ has order $3$ if 
$i=1,2,\do,n-2$ and order $4$ if $i=n-1$; we have $s_is_j=s_js_i$ if $|i-j|\ge2$.

(i) {\it We have $\ufK^J=0$ unless}

($*$) {\it $J=\{s_{k+1},s_{k+2},\do,s_n\}$ with $0\le k\le n$ such that $n-k=a^2+a$ for some 
$a\in\NN$.}
\nl
Now assume that $(*)$ holds. Then $\cw$ is a Weyl group of type $B_k$ with standard generators
$\s_1,\s_2,\do,\s_k$ where $\s_i=s_i$ for $1\le i<k$ and $\s_k$ is the unique element in 
the subgroup of $\WW$ generated by $s_k,s_{k+1},\do,s_n$ such that $\s_k\in\cw-\{1\}$ (if 
$k\ge1$). Let $\tl:\cw@>>>\NN$ be the length function of the Weyl group $\cw$.

(ii) {\it For any $w\in\cw$ we have $\dim({}^w\ufK^J)=1$.}

(iii) {\it For any $i\in[1,k]$ there is a unique element $x\in{}^{\s_i}\ufK^J-\{0\}$ such that
$(x+\bold1)*(x-c\bold1)=0$ for some $c\in v\ZZ[v]$; in fact we have $c=v^2$ if $1\le i<k$ 
and $c=v^{4a+2}$ if $i=k$. We set $x=t_i$.}

(iv) {\it For any $w\in\cw$ there is a unique element $t_w\in{}^w\ufK^J-\{0\}$ such that the
following hold: $t_{\s_i}=t_i$ for $i\in[1,k]$; $t_w*t_{w'}=t_{w'w}$ if $w,w'\in\cw$,
$\tl(w'w)=\tl(w)+\tl(w')$.}

We see that $\ufK^J$ is an Iwahori-Hecke algebra with not necessarily equal
parameters. Similar results hold for other classical groups.

\head 37. A Mackey type formula\endhead  
\subhead 37.1\endsubhead
We fix a connected component $D$ of $G$. With notation in 26.1, if 
$J,J'\sub\II,P\in\cp_J,Q\in\cp_{J'},u=\po(P,Q)\in\WW$, then $u\in{}^J\WW^{J'}$. Setting 
$P^Q=(P\cap Q)U_P$, we have $P^Q\in\cp_{J\cap\Ad(u)J'}$. 

For $K,K'\sub\II$ and $u\in{}^K\WW^{K'}$ let  
$$\align\Up_u=&\{(X,Y,g(U_X\cap U_Y));X\in\cp_{K\cap\Ad(u)K'},
Y\in\cp_{K'\cap\Ad(u\i)K},\\&g(U_X\cap U_Y)\in D/(U_X\cap U_Y),\po(X,Y)=u\}.\endalign$$
We have a diagram
$$Z_{K\cap\Ad(u)K',D}@<\fj<<\Up_u@>\fh>>Z_{K'\cap\Ad(u\i)K,D}$$
where $\fj(X,Y,g(U_X\cap U_Y))=(X,gU_X),\fh(X,Y,g(U_X\cap U_Y))=(Y,gU_Y)$. Set
$$\Ph_u=\fh_!\fj^*:\cd(Z_{K\cap\Ad(u)K',D})@>>>\cd(Z_{K'\cap\Ad(u\i)K,D}).$$

\proclaim{Proposition 37.2}Let $K,K',J\sub\II$ be such that $K\sub J,K'\sub J$. Let 
$A'\in\cd(Z_{K,D})$. We set $\fB=\fe_{K',J}\ff_{K,J}A'\in\cd(Z_{K',D})$. For any 
$u\in{}^K\WW^{K'}\cap\WW_J$ we set 
$\fB_u=\ff_{K'\cap\Ad(u\i)K,K'}\Ph_u\fe_{K\cap\Ad(u)K',K}A'\in\cd(Z_{K',D})$ and
$m_u=\dim(U_P\cap U_R)/U_Q$ where $P\in\cp_K,R\in\cp_{K'}$, $\po(P,R)=u$ and
$Q=Q_{J,P}=Q_{J,R}\in\cp_J$ (notation of 36.4). We have 
$$\fB\Bpq\{\fB_u[[-m_u]];u\in{}^K\WW^{K'}\cap\WW_J\},$$
with $\Bpq$ as in 32.15.
\endproclaim
We have a commutative diagram with a cartesian square
$$\CD 
{} @.   \fE@>\fb>>    Z_{K',J,D}@>\fc'>>Z_{K',D}\\
@.      @V\fa VV       @V\fd'VV         @.   \\
Z_{K,D}@<\fc<<Z_{K,J,D}@>\fd>>Z_{J,D}@.         {}                               \endCD$$
Here 

$\fE=\{(P,R,gU_Q);P\in\cp_K,R\in\cp_{K'},gU_Q\in D/U_Q,Q=Q_{J,P}=Q_{J,R}\}$,

$\fc(P,gU_Q)=(P,gU_P)$, $\fd(P,gU_Q)=(Q,gU_Q)$ with $Q=Q_{J,P}$,

$\fc'(R,gU_Q)=(R,gU_R)$, $\fd'(R,gU_Q)=(Q,gU_Q)$ with $Q=Q_{J,R}$,

$\fa(P,R,gU_Q)=(P,gU_Q),\fb(P,R,gU_Q)=(R,gU_Q)$.
\nl
We have 
$$\fB=\fc'_!\fd'{}^*\fd_!\fc^*A'=\fc'_!\fb_!\fa^*\fc^*A'=(\fc'\fb)_!(\fc\fa)^*A'
=\fq_!\fp^*A'$$
where $\fq=\fc'\fb:\fE@>>>Z_{K',D},\fp=\fc\fa:\fE@>>>Z_{K,D}$ are given by 
$\fq(P,R,gU_Q)=(R,gU_R)$, $\fp(P,R,gU_Q)=(P,gU_P)$. We have a partition
$$\fE=\sqc_{u\in{}^K\WW^{K'}\cap\WW_J}\fE_u$$ 
where $\fE_u=\{(P,R,gU_Q)\in\fE;\po(P,R)=u\}$ is locally closed in $\fE$. Let 
$\fp_u=\fp|_{\fE_u}:\fE_u@>>>Z_{K,D},\fq_u=\fq|_{\fE_u}:\fE_u@>>>Z_{K'}$. By 32.15, we have
$$\fq_!\fp^*A'\Bpq\{\fq_{u!}\fp_u^*A';u\in{}^K\WW^{K'}\cap\WW_J\}.$$
It remains to show that, for $u$ as above, we have
$$\fq_{u!}\fp_u^*A'=\fB_u[[-m_u]].$$
We have a commutative diagram
$$\CD 
Z_{K,D}@<\ti\fd<<Z_{K\cap\Ad(u)K',K,D}@>\ti\fc>>Z_{K\cap\Ad(u)K',D} @.   {}   \\
  @.              @A\ti\fa AA                   @A\fj AA                 @.   \\
  {}     @.      \fS'_u    @>\ti\fb>>            \Up_u@>\fh>>Z_{K'\cap\Ad(u\i)K,D}  \\
  @.          @A\fx AA                 @.           @A\ti\fc'AA    \\
  {} @.   \fS_u     @>\fy>>        {}  @>>>  Z_{K'\cap\Ad(u\i)K,K',D}\\
  @.       @.                      @.          @V\ti\fd'VV               \\
  {} @.    {}@.                    {}@.        Z_{K',D}                  \endCD$$
Here 

$\ti\fc(X,gU_P)=(X,gU_X),\ti\fd(X,gU_P)=(P,gU_P)$,

$\ti\fc'(Y,gU_R)=(Y,gU_Y),\ti\fd'(Y,gU_R)=(R,gU_R)$,

$\fS'_u,\ti\fa,\ti\fb$ are defined so that the square $(\ti\fa,\ti\fb,\ti\fc,\fj)$ is 
cartesian; 

$\fS_u,\fx,\fy$ are defined so that the square $(\fx,\fy,\fh\ti\fb,\ti\fc')$ is cartesian.
\nl
Then 
$$\align\fS_u=&\{(X,Y,g(U_X\cap U_Y),g'U_P,g''U_R);X\in\cp_{K\cap\Ad(u)K'},
Y\in\cp_{K'\cap\Ad(u\i)K},\\&g(U_X\cap U_Y)\in D/(U_X\cap U_Y),g'U_P\in D/U_P,
g''U_R\in D/U_R,\po(X,Y)=u\\&,P=Q_{K,X},R=Q_{K',Y},g'U_X=gU_X,g''U_Y=gU_Y\}.\endalign$$
Set $\fr=\ti\fd'\fy:\fS_u@>>>Z_{K',D},\fs=\ti\fd\ti\fa\fx:\fS_u@>>>Z_{K,D}$. Then
$$\align&\fr(X,Y,g(U_X\cap U_Y),g'U_P,g''U_R)=(P,g'U_P),\\&\fs(X,Y,g(U_X\cap U_Y),
g'U_P,g''U_R)=(R,g''U_R).\endalign$$
We have 
$$\align&\fB_u=\ti\fd'_!\ti\fc'{}^*\fh_!\fj^*\ti\fc_!\ti\fd^*A'
=\ti\fd'_!\ti\fc'{}^*\fh_!\ti\fb_!\ti\fa^*\ti\fd^*A'\\&=
\ti\fd'_!\ti\fc'{}^*(\fh\ti\fb)_!(\ti\fd\ti\fa)^*A'
=\ti\fd'_!\fy_!\fx^*(\ti\fd\ti\fa)^*A'=(\ti\fd'\fy)_!(\ti\fd\ti\fa\fx)^*A'=\fr_!\fs^*A'.
\endalign$$
We show that $\fq_{u!}\fp_u^*A'=\fr_!\fs^*A'[[-m_u]]$. We have a commutative diagram
$$\CD 
Z_{K,D}@<\fp_u<<\fE_u@>\fq_u>>Z_{K',D}\\
@V1VV   @V\ft VV    @V1VV\\
Z_{K,D}@<\fs<<\fS_u@>\fr>>Z_{K',D}           \endCD$$
where $\ft(P,R,gU_Q)=(P^R,R^P,g(U_{P^R}\cap U_{R^P}),gU_P,gU_R)$ is well defined since 
$U_Q\sub U_P,U_Q\sub U_R,U_Q\sub U_{P^R},U_Q\sub U_{R^P}$. We continue the proof assuming
that

(a) {\it $\ft$ is an affine space bundle with fibres of dimension $m_u$.}
\nl
For any $\tA\in\cd(\fS_u)$ we have $\ft_!\ft^*(\tA)=\tA[[-m_u]]$. Hence 
$$\fr_!\fs^*A'[[-m_u]]=\fr_!\ft_!\ft^*\fs^*A'=(\fr\ft)_!(\fs\ft)^*A'=\fq_{u!}\fp_u^*A',$$
as required. 

We prove (a). We only show that each fibre of $\ft$ is an affine space of dimension $m_u$.
First we show that $\ft:\fE_u@>>>\fS_u$ is surjective. Assume that we are given
$(X,Y,g_0(U_X\cap U_Y),g'_0U_P,g''_0U_R)\in\fS_u$. Then $g_0,g'_0,g''_0\in D$, 
$v'=g_0\i g'_0\in U_X$, $v''=g_0\i g''_0\in U_Y$. We must show that there exists $g\in D$
such that $g_0\in g(U_X\cap U_Y),g'_0\in gU_P$, $g''_0\in gU_R$. Setting $y=g\i g_0$, we 
must show that there exists $y\in U_X\cap U_Y$ such that $yv'\in U_P,yv''\in U_R$. We have
$v'=v'_1v'_2$ where $v'_1\in U_R\cap P$, $v'_2\in U_P$ and $v''=v''_1v''_2$ where 
$v''_1\in U_P\cap R$, $v''_2\in U_R$. Then $v'_1\in U_X\cap U_Y,v''_1\in U_X\cap U_Y$. 
Setting $y=(v'_1v''_1)\i\in U_X\cap U_Y$ we have $yv'=v''_1{}\i v'_2\in U_P$ and 
$yv''=v''_1{}\i v'_1{}\i v''_1v''_2\in v''_1{}\i U_Rv''_1U_R=U_R$, as desired. 

It remains to show that, if $(P,R,gU_Q)\in\fE_u$, then

(b) {\it $F=\{(P',R',g'U_{Q'})\in E_u;\ft(P',R',g'U_{Q'})=\ft(P,R,gU_Q)\}$ is an affine 
space of dimension $m_u$.}
\nl
For $(P',R',g'U_{Q'})\in F$, both $P,P'$ contain $P^R=(P')^{R'}$ and have the same type
hence $P=P'$. Similarly $R=R',Q=Q'$. Hence 
$$\align&
F\cong\{g'U_Q;gU_P=g'U_P,gU_R=g'U_R,g(U_{P^R}\cap U_{R^P})=g'(U_{P^R}\cap U_{R^P})\}\\&
=\{g'U_Q;g\i g'\in U_P\cap U_R\cap U_{P^R}\cap U_{R^P}\}\\&
=\{g'U_Q;g\i g'\in U_P\cap U_R\}\cong(U_P\cap U_R)/U_Q,\endalign$$
and (b) follows. This completes the proof.

\subhead 37.3\endsubhead
In the remainder of this section we assume that we are in the setup of 36.8. From 37.2 we 
deduce that for $K,K',J\sub\II$ such that $K\sub J,K'\sub J$ we have
$$\fe_{K',J}\ff_{K,J}=\sum_{u\in{}^K\WW^{K'}\cap\WW_J}
v^{m_u}\ff_{K'\cap\Ad(u\i)K,K'}\Ph_u\fe_{K\cap\Ad(u)K',K}$$
as $\ca$-linear maps $\fK(Z_{K,D})@>>>\fK(Z_{K',D})$, where the $\ca$-linear map 

$\fK(Z_{K\cap\Ad(u)K',D})@>>>\fK(Z_{K'\cap\Ad(u\i)K,D})$
\nl
induced by $\Ph_u$ in 37.1 is denoted again by $\Ph_u$. 

\subhead 37.4\endsubhead
Let $K,K'\sub J$ and let $x\in\fK(Z_{K,D})^K,x'\in\fK(Z_{K',D})^{K'}$. We show: 
$$(\ff_{K,J}x:\ff_{K',J}x')=\sum_{u\in{}^K\WW^{K'}\cap\WW_J;K=\Ad(u)K'}v^{m(u)}(\Ph_ux:x').
\tag a$$
Here $m(u),\Ph_u$ are as in 37.3; in our case, $\Ph_u:\fK(Z_{K,D})@>>>\fK(Z_{K',D})$. 
Using 37.3 we have
$$\align&(\ff_{K,J}x:\ff_{K',J}x')=(\fe_{K',J}\ff_{K,J}x:x')\\&=
(\sum_{u\in{}^K\WW^{K'}\cap\WW_J}v^{m_u}\ff_{K'\cap\Ad(u\i)K,K'}\Ph_u
\fe_{K\cap\Ad(u)K',K}x:x')\\&=\sum_{u\in{}^K\WW^{K'}\cap\WW_J}
v^{m_u}(\Ph_u\fe_{K\cap\Ad(u)K',K}x:\fe_{K'\cap\Ad(u\i)K,K'}x')\\&
=\sum_{u\in{}^K\WW^{K'}\cap\WW_J;K\cap\Ad(u)K'=K,K'\cap\Ad(u\i)K=K'}v^{m_u}(\Ph_ux:x').
\endalign$$
The condition that $K\cap\Ad(u)K'=K,K'\cap\Ad(u\i)K=K'$ is equivalent to $K\sub\Ad(u)K'$,
$K'\sub\Ad(u\i)K$ that is, to $K=\Ad(u)K'$. This proves (a).

Consider the equivalence relation on the set of subsets of $J$ given by $K_1\si K_2$ if 
$\Ad(u)K_1=K_2$ for some $u\in\WW_J$. For any equivalence class $\fo$ under $\si$ we set 
$$\ufK(Z_{J,D})^\fo=\sum_{H\sub J;H\in\fo}f_{H,J}(\ufK(Z_{H,D})^H)\sub\ufK(Z_{J,D}).$$ 
Note that $\ufK(Z_{J,D})^\fo$ is homogeneous. 

\proclaim{Proposition 37.5} We have
$$\ufK(Z_{J,D})=\op_\fo\ufK(Z_{J,D})^\fo\tag a$$
where $\fo$ runs over the equivalence classes for $\si$.
\endproclaim
If $J=\em$ we have $\ufK(Z_{J,D})=\ufK(Z_{J,D})^J$ and the result is obvious. We may assume
that $J\ne\em$ and that the result is true when $J$ is replaced by a strictly smaller 
subset. Using 36.13(d) and the induction hypothesis we have
$$\align&\ufK(Z_{J,D})=\ufK(Z_{J,D})^J+\sum_{L\subsetneq J}\ff_{L,J}\sum_{J';J'\sub L}
\ff_{J',L}\ufK(Z_{J',D})^{J'}\\&\sub\sum_{J';J'\sub J}\ff_{J',J}\ufK(Z_{J',D})^{J'}.
\endalign$$
Thus, $\ufK(Z_{J,D})=\sum_\fo\ufK(Z_{J,D})^\fo$. Next we show that
$(\ufK(Z_{J,D})^\fo:\ufK(Z_{J,D})^{\fo'})=0$ if $\fo\ne\fo'$. It is enough to show that 
$(\ff_{H,J}(\fK(Z_{H,D})^H):\ff_{H',J}(\fK(Z_{H',D})^{H'}))=0$ if $H,H'\sub J$, 
$H\not\si H'$. This follows from 37.4(a). It remains to use the following (easily verified)
statement: if $V$ is a finite dimensional vector space with a non-singular symmetric 
bilinear form $(:)$ and $V_1,V_2,\do,V_k$ are subspaces such that $(V_i:V_j)=0$ for 
$i\ne j$ and $\sum_iV_i=V$ then $V=\op_iV_i$. 

\head 38. Duality\endhead
\subhead 38.1\endsubhead
We fix a connected component $D$ of $G$ that generates $G$. In this section we study an 
involution of the set of isomorphism classes of character sheaves on $D$ called {\it 
duality}. 

We write $\e$ instead of $\e_D:\WW@>>>\WW$ (see 26.2). For $J\sub\II$ such that $\e(J)=J$ 
we set, as in 30.3: 
$$V_{J,D}=\{(P,gU_P);P\in\cp_J,gU_P\in N_DP/U_P\}.$$
We also set $\WW_J^\e=\{w\in\WW_J;\e(w)=w\}$ where $\WW_J$ is as in 26.1. Let $J_\e$ be the
set of orbits of the restriction of $\e$ to $J$.

For any $P_0\in\cp_J$ we have a functor $A_0\m A_0^\flat$ (see 30.3) from the category of
$P_0/U_{P_0}$-equivariant perverse sheaves on the connected component $N_DP_0/U_{P_0}$ of
$N_GP_0/U_{P_0}$ to the category of perverse sheaves on $V_{J,D}$. 

Let $CS(V_{J,D})$ be the full subcategory of the category of perverse sheaves on 
$V_{J,D}$ whose objects are isomorphic to objects of the form $A_0^\flat$ where $A_0$ is a
direct sum of character sheaves on $N_DP_0/U_{P_0}$. In particular, $CS(V_{\II,D})$ 
(denoted also by $CS(D)$) is the category of perverse sheaves on $D$ that are direct sums
of character sheaves. Note that $A_0\m A_0^\flat$ is an equivalence of categories 
$CS(N_DP_0/U_{P_0})@>\si>>CS(V_{J,D})$.

For $J\sub J'\sub\II$ such that $\e(J)=J,\e(J')=J'$ we have functors
$f_{J,J'}:\cd(V_{J,D})@>>>\cd(V_{J',D})$ and $e_{J,J'}:\cd(V_{J',D})@>>>\cd(V_{J,D})$, see
30.4. From definitions, for $J\sub J'\sub J''\sub\II$ such that 
$\e(J)=J,\e(J')=J',\e(J'')=J''$, we have
$$f_{J,J''}=f_{J',J''}f_{J,J'},\qua e_{J,J''}=e_{J,J'}e_{J',J''}.\tag a$$
Clearly, $f_{J,J}=1,e_{J,J}=1$.

\subhead 38.2\endsubhead
We show that:

(a) {\it for $J\sub J'$ as above, $e_{J,J'}$ restricts to a functor 
$CS(V_{J',D})@>>>CS(V_{J,D})$ denoted again by $e_{J,J'}$.}
\nl
Let $P\in\cp_J,P'\in\cp_{J'}$ be such that $P\sub P'$. Let $D_0=N_DP/U_P$,
$D'_0=N_DP'/U_{P'}$. Let $C_0\in CS(D'_0)$ and let $C=C_0^\flat$ be the corresponding 
object of $CS(V_{J',D})$. From 31.14 we see that $\res_{D'_0}^{D_0}(C_0)\in CS(D_0)$ and 
from 30.4(b) we see that $e_{J,J'}C$ is the perverse sheaf 
$(\res_{D'_0}^{D_0}(C_0))^\flat$ on $V_{J,D}$ hence $e_{J,J'}C\in CS(V_{J,D})$. This proves
(a).

\subhead 38.3\endsubhead
We show that:

(a) {\it for $J\sub J'$ as above, $f_{J,J'}$ restricts to a functor 
$CS(V_{J,D})@>>>CS(V_{J',D})$ denoted again by $f_{J,J'}$.}
\nl
Let $P,P',D_0,D'_0$ be as in 38.2. Let $A_0\in CS(D_0)$ and let $A=A_0^\flat$ be the 
corresponding object of $CS(V_{J,D})$. By 30.4(a), $f_{J,J'}A=A'_0{}^\flat$ where 
$A'_0=\ind_{D_0}^{D'_0}A_0$ is a direct sum of simple admissible perverse sheaves on 
$D'_0$. It remains to show that $A'_0\in CS(D'_0)$. To do this we may assume that 
$J'=\II$ hence $P'=G^0,D'_0=D$. Let $\a=\dim U_P$. Let $L$ be a Levi of $P$. We can 
identify naturally $N_GP/U_P$ with $H=N_GP\cap N_GL$, a reductive group with $H^0=L$. Then
$D_0$ becomes $H\cap D$. We identify the canonical torus of $L$ with the canonical torus 
$\TT$ of $G^0$, and the Weyl group of $L$ with the subgroup $\WW_J$ of $\WW$ as in 29.1. 
Let $\cl\in\fs(\TT)$. Let $\ss=(s_1,s_2,\do,s_r)$ be a sequence in $J\cup\{1\}$ such that 
$s_1s_2\do s_r\uD\in\Wb_\cl$ (see 28.3). Let $\bZ^\ss_{\em,\II,D}$, 
$\bK^{\ss,\cl}_D\in\cd(D)$ be as in 28.12, $\bcl$ as in 28.9, and let 
$\bZ^\ss_{\em,J,D_0}$, $\bK^{\ss,\cl}_{D_0}\in\cd(D_0)$, $\bcl_0$, be the analogous
objects defined in terms of $H,D_0$ instead of $G,D$. Consider the commutative diagram
$$\CD
\bZ^\ss_{\em,J,D_0}@<pr_2<<V_1\T_{D_0}\bZ^\ss_{\em,J,D_0}@>f_0>>\bZ^\ss_{\em,\II,D}\\
@Vf_1VV                                      @Vpr_1VV      @Vf_2VV \\
D_0@<a_1<<                                  V_1@>a'>>      V_2@>a''>>D          \endCD$$
where

$V_1=\{(g,x)\in D\T G^0;x\i gx\in N_DP\}$,

$V_2=\{(g,xP)\in D\T G^0/P;x\i gx\in N_DP\}$,

$a_1(g,x)=g'$ where $g_0\in D_0$ is such that $x\i gx\in g_0U_P$,

$a'(g,x)=(g,xP),a''(g,xP)=g$, $f_1(\b_0,\b_1,\do,\b_r,g_0)=g_0$,

$f_2(B_0,B_1,\do,B_r,g)=(g,x_0P)$ where $x_0\in G^0$ is such that $x_0\i B_0x_0\sub P$,

$f_0((g,x),(\b_0,\b_1,\do,\b_r,g_0))=(x\b_0U_Px\i,x\b_1U_Px\i,\do,x\b_rU_Px\i,g)$.
\nl
Both squares in the diagram are cartesian and the maps $a',f_0$ (resp. $a_1,pr_2$) are 
smooth with connected fibres of dimension $\dim G-\a$ (resp. $\dim G+\a$). From definitions
we have $pr_2^*\bcl'=f_0^*\bcl$. Hence
$$a_1^*(\bK^{\ss,\cl}_{D_0})=
     a_1^*f_{1!}\bcl_0= pr_{1!}pr_2^*\bcl_0=pr_{1!}f_0^*\bcl= a'{}^*f_{2!}\bcl.$$
We see that 
$$a_1^\bst(\bK^{\ss,\cl}_{D_0})[-\dim G-\a]=a'{}^\bst\tK[-\dim G+\a]$$ 
where $\tK=f_{2!}\bcl$, that is $a_1^\bst(\bK^{\ss,\cl}_{D_0})=a'{}^\bst\tK[2\a]$. Hence
$$\align&a_1^\bst({}^pH^i(\bK^{\ss,\cl}_{D_0}))={}^pH^i(a_1^\bst(\bK^{\ss,\cl}_{D_0}))
={}^pH^i(a'{}^\bst\tK[2\a])=a'{}^\bst({}^pH^i(\tK[2\a]))\\&=a'{}^\bst({}^pH^{i+2\a}\tK).
\endalign$$
From this and definition (27.1) we see that 
$$\ind_{D_0}^D({}^pH^i(\bK^{\ss,\cl}_{D_0}))=a''_!({}^pH^{i+2\a}\tK).$$
We have 
$$\op_i\ind_{D_0}^D({}^pH^i(\bK^{\ss,\cl}_{D_0}))[-i]=
\op_i{}^pH^{i+2\a}(\bK^{\ss,\cl}_D)[-i]\in\cd(D).\tag b$$
Indeed, the left hand side is
$$\align&\op_ia''_!({}^pH^{i+2\a}\tK)[-i]=a''_!(\tK[2\a])=a''_!f_{2!}\bcl[2\a]=
\bK^{\ss,\cl}_D)[2\a]\\&=\op_i{}^pH^{i+2\a}(\bK^{\ss,\cl}_D)[-i]\in\cd(D),\endalign$$
where we use that $\tK$ and $\bK^{\ss,\cl}_D$ are semisimple complexes (a consequence of
the decomposition theorem \cite{\BBD}). Since ${}^pH^i(\bK^{\ss,\cl}_{D_0})$ is a direct 
sum of character sheaves on $D_0$, we see, using 30.6(a), that 
$\ind_{D_0}^D({}^pH^i(\bK^{\ss,\cl}_{D_0}))$ is a perverse sheaf on $D$. Taking ${}^pH^i$
for both sides of (b) we therefore find
$$\ind_{D_0}^D({}^pH^i(\bK^{\ss,\cl}_{D_0}))={}^pH^{i+2\a}(\bK^{\ss,\cl}_D)\tag c$$
for any $i\in\ZZ$. To prove (a) it is enough to verify the following statement.

(d) {\it If $A_1\in\hD_0^\cl$, then $\ind_{D_0}^D(A_1)$ is a direct sum of character
sheaves in $\hD^\cl$.}
\nl
We may assume that $A_1$ is a direct summand of ${}^pH^i(\bK^{\ss,\cl}_{D_0})$. Then 
$\ind_{D_0}^D(A_1)$ is a direct summand of $\ind_{D_0}^D({}^pH^i(\bK^{\ss,\cl}_{D_0}))$. 
From (c) we see that $\ind_{D_0}^D(A_1)$ is a direct summand of 
${}^pH^{i+2\a}(\bK^{\ss,\cl}_D)$ which is a direct sum of character sheaves in $\hD^\cl$.
This proves (d) hence also (a).

\subhead 38.4\endsubhead
For $J\sub J'\sub\II$ such that $\e(J)=J,\e(J')=J'$ we have functors
$\tf_{J,J'}:\cd(V_{J,D})@>>>\cd(V_{J',D})$, $\te_{J,J'}:\cd(V_{J',D})@>>>\cd(V_{J,D})$, see
30.4. 

Let $K,K',J\sub\II$ be such that $K\sub J,K'\sub J,\e(J)=J,\e(K)=K,\e(K')=K'$. For any 
$u\in{}^K\WW^{K'}\in\WW_J^\e$ let 
$$\align\Xi_u=&\{(X,Y,g(U_X\cap U_Y));X\in\cp_{K\cap\Ad(u)K'},Y\in\cp_{K'\cap\Ad(u\i)K},
\\&g(U_X\cap U_Y)\in(N_DX\cap N_DY)/(U_X\cap U_Y),\po(X,Y)=u\}.\endalign$$ 
We have a diagram
$$V_{K\cap\Ad(u)K',D}@<j<<\Xi_u@>h>>V_{K'\cap\Ad(u\i)K,D}$$
where $j(X,Y,g(U_X\cap U_Y))=(X,gU_X),h(X,Y,g(U_X\cap U_Y))=(Y,gU_Y)$. Set
$\Ps_u=h_!j^*:\cd(V_{K\cap\Ad(u)K',D})@>>>\cd(V_{K'\cap\Ad(u\i)K,D})$.

\proclaim{Lemma 38.5}Let $A'\in\cd(V_{K,D})$. We set 
$\ti{\fC}=\te_{K',J}\tf_{K,J}A'\in\cd(V_{K',D})$. For any 
$u\in{}^K\WW^{K'}\cap\WW_J^\e$ we set 
$\ti{\fC}_u=\tf_{K'\cap\Ad(u\i)K,K'}\Ps_u\te_{K\cap\Ad(u)K',K}A'\in\cd(V_{K',D})$ and
$m_u=\dim(U_P\cap U_R)/U_Q$ where $P\in\cp_K,R\in\cp_{K'}$, $\po(P,R)=u$ and 
$Q=Q_{J,P}=Q_{J,R}\in\cp_J$ (notation of 36.4). We have 
$$\ti{\fC}\Bpq\{\ti{\fC}_u[[-m_u]];u\in{}^K\WW^{K'}\cap\WW_J^\e\}$$
where $\Bpq$ is as in 32.15. 
\endproclaim
The proof is very similar to that of Proposition 37.2. We have a commutative diagram with a
cartesian square
$$\CD 
{} @.    E@>b>>    V_{K',J,D}@>c'>>V_{K',D}\\
@.      @VaVV       @Vd'VV         @.   \\
V_{K,D}@<c<<V_{K,J,D}@>d>>V_{J,D}@.         {}                       \endCD$$
Here 

$E=\{(P,R,gU_Q);P\in\cp_K,R\in\cp_{K'},gU_Q\in(N_DP\cap N_DR)/U_Q,Q=Q_{J,P}=Q_{J,R}\}$,

$c(P,gU_Q)=(P,gU_P)$, $d(P,gU_Q)=(Q,gU_Q)$ with $Q=Q_{J,P}$,

$c'(R,gU_Q)=(R,gU_R)$, $d'(R,gU_Q)=(Q,gU_Q)$, with $Q=Q_{J,R}$,

$a(P,R,gU_Q)=(P,gU_Q),b(P,R,gU_Q)=(R,gU_Q)$.
\nl
We have 
$$\ti{\fC}=c'_!d'{}^*d_!c^*A'=c'_!b_!a^*c^*A'=(c'b)_!(ca)^*A'=q_!p^*A'$$
where $q=c'b:E@>>>V_{K',D},p=ca:E@>>>V_{K,D}$ are given by $q(P,R,gU_Q)=(R,gU_R)$, 
$p(P,R,gU_Q)=(P,gU_P)$. We have a partition
$$E=\sqc_{u\in{}^K\WW^{K'}\cap\WW_J^\e}E_u$$ 
where $E_u=\{(P,R,gU_Q)\in E;\po(P,R)=u\}$ is locally closed in $E$. Let 
$p_u=p|_{E_u}:E_u@>>>V_{K,D},q_u=q|_{E_u}:E_u@>>>V_{K',D}$. By 32.15, we have
$$q_!p^*A'\Bpq\{q_{u!}p_u^*A';u\in{}^K\WW^{K'}\cap\WW_J^\e\}.$$
It remains to show that, for $u$ as above, we have $q_{u!}p_u^*A'=\ti{\fC}_u[[-m_u]]$. We
have a commutative diagram
$$\CD 
V_{K,D}@<\td<<V_{K\cap\Ad(u)K',K,D}@>\tc>>V_{K\cap\Ad(u)K',D} @.   {}      \\
  @.          @A\ta AA                @AjAA               @.            \\
  {}     @.   S'_u    @>\tb>>       \Xi_u@>h>>V_{K'\cap\Ad(u\i)K,D}  \\
  @.          @AxAA                 @.                   @A\tc'AA    \\
  {} @.   S_u     @>y>>        {}  @>>>   V_{K'\cap\Ad(u\i)K,K',D}\\
  @.      @.                   @.               @V\td'VV       \\
  {}  @.  {}  @.               {} @.           V_{K',D} \endCD$$
Here 

$\tc(X,gU_P)=(X,gU_X),\td(X,gU_P)=(P,gU_P)$,

$\tc'(Y,gU_R)=(Y,gU_Y),\td'(Y,gU_R)=(R,gU_R)$,

$S'_u,\ta,\tb$ are defined so that the square $(\ta,\tb,\tc,j)$ is cartesian; 

$S_u,x,y$ are defined so that the square $(x,y,h\tb,\tc')$ is cartesian.
\nl
Then 
$$\align&S_u=\{(X,Y,g(U_X\cap U_Y),g'U_P,g''U_R);X\in\cp_{K\cap\Ad(u)K'},
Y\in\cp_{K'\cap\Ad(u\i)K},\\&g(U_X\cap U_Y)\in(N_DX\cap N_DY)/(U_X\cap U_Y),
g'U_P\in N_DX/U_P,\\&g''U_R\in N_DY/U_R,\po(X,Y)=u,P=Q_{K,X},R=Q_{K',Y},\\&
g'U_X=gU_X,g''U_Y=gU_Y\}.\endalign$$
Set $r=\td'y:S_u@>>>V_{K',D},s:\td\ta x:S_u@>>>V_{K,D}$. Then
$$\align&r(X,Y,g(U_X\cap U_Y),g'U_P,g''U_R)=(P,g'U_P),\\&
s(X,Y,g(U_X\cap U_Y),g'U_P,g''U_R)=(R,g''U_R).\endalign$$
We have 
$$\align&\ti{\fC}_u=\td'_!\tc'{}^*h_!j^*\tc_!\td^*A'
=\td'_!\tc'{}^*h_!\tb_!\ta^*\td^*A'=\td'_!\tc'{}^*(h\tb)_!(\td\ta)^*A'\\&
=\td'_!y_!x^*(\td\ta)^*A'=(\td'y)_!(\td\ta x)^*A'=r_!s^*A'.\endalign$$
We show that $q_{u!}p_u^*A'=r_!s^*A'[[-m_u]]$. We have a commutative diagram
$$\CD 
V_{K,D}@<p_u<<E_u@>q_u>>V_{K',D}\\
@V1VV   @Vt VV    @V1VV\\
V_{K,D}@<s<<S_u@>r>>V_{K',D}                                 \endCD$$
where $t(P,R,gU_Q)=(P^R,R^P,g(U_{P^R}\cap U_{R^P}),gU_P,gU_R)$ is well defined since 
$U_Q\sub U_P,U_Q\sub U_R,U_Q\sub U_{P^R},U_Q\sub U_{R^P}$. We continue the proof assuming
that

(a) {\it $t$ is an affine space bundle with fibres of dimension $m_u$.}
\nl
For any $\tA\in\cd(S_u)$ we have $t_!t^*(\tA)=\tA[[-m_u]]$. Hence 
$$r_!s^*A'[[-m_u]]=r_!t_!t^*s^*A'=(rt)_!(st)^*A'=q_{u!}p_u^*A',$$
as required. 

We prove (a). Consider the commutative diagram
$$\CD      
E_u@>>>\fE_u\\
@VtVV @V\ft VV\\
S_u@>>>\fS_u                                  \endCD$$
where $\fE_u@>\ft>>\fS_u$ is as in the proof of 37.2 and the horizontal maps are the 
obvious imbeddings. Clearly, this diagram is cartesian. Hence (a) is a consequence of the
analogous statement 37.2(a) for $\fE_u@>\ft>>\fS_u$. This completes the proof.

\subhead 38.6\endsubhead
In the setup of 38.4, let $u\in{}^K\WW^{K'}\cap\WW_J^\e$. Let 
$\d=\dim(U_{X_u}/(U_{X_u}\cap U_{Y_u}))=\dim(U_{Y_u}/(U_{X_u}\cap U_{Y_u}))$ where 
$(X_u,Y_u,g(U_{X_u}\cap U_{Y_u}))\in\Xi_u$. Let $\a_u=\dim U_{X_u}=\dim U_{Y_u}$, 
$\a=\dim U_P,\a'=\dim U_R$ where $P\in\cp_K,R\in\cp_{K'}$. We show:

(a) {\it $h$ and $j$ in 38.4 are affine space bundles with fibres of dimension $\d$.}
\nl
It is enough to prove the statements on $j$ (the statement on $h$ is entirely similar). We
show that $j$ is surjective. Let $X,g$ be such that $(X,gU_X)\in V_{K\cap\Ad(u)K',D}$. We 
must show that there exist $Y\in\cp_{K'\cap\Ad(u\i)K},g'\in N_DX\cap N_DY$ such that 
$\po(X,Y)=u$ and $g'U_X=gU_X$. Setting $g\i g'=v$, it is enough to show that for any 
$Y\in\cp_{K'\cap\Ad(u\i)K}$ with $\po(X,Y)=u$ there exists $v\in U_X$ such that
$gv\in N_GY$. Now $X,Y$ contain a common Levi $M$. Since $g\in N_GX$ we can find $v\in U_X$
such that $g'=gv\in N_GX\cap N_GM$. There is a unique parabolic $Y'$ of the same type as 
$Y$ such that $Y'$ has Levi $M$ and $\po(X,Y')=u$. Then $\po(g'Xg'{}\i,g'Yg'{}\i)=u$, 
$g'Yg'{}\i$ has Levi $g'Mg'{}\i=M$. By uniqueness, we have $Y'=Y$. Thus $g'\in N_GY$.

We show that the fibres of $j$ are affine spaces of dimension $\d$. Let 
$(X,Y,g(U_X\cap U_Y))\in\Xi_u$. We must show that 
$F=\{(X,Y',g'(U_X\cap U_Y'))\in\Xi_u;gU_X=g'U_X\}$ is an affine space. Fix 
$Y'\in\cp_{K'\cap\Ad(u\i)K}$ such that $\po(X,Y')=u$. (The set of such $Y'$ is a 
homogeneous space $U_X/(U_X\cap Y)$, hence is an affine space of dimension $\d$.) It is
enough to show that

$\{g'(U_X\cap U_Y')\in(N_DX\cap N_DY)/(U_X\cap U_Y);g'\in gU_X\}$ is a point.
\nl
Now $g=g_0v_0$ where $g_0\in N_DX\cap N_DY,v_0\in U_X$. It is enough to show that 
$\{v\in U_X;v_0v\in N_GX\cap N_GY\}/(U_X\cap U_Y)$ is a point or that
$(U_X\cap N_GX\cap N_GY)/(U_X\cap U_Y)$ is a point or that $U_X\cap Y=U_X\cap U_Y$. This is
clear.

\subhead 38.7\endsubhead
Let $u\in{}^K\WW^{K'}\cap\WW_J^\e$. We set $H=K\cap\Ad(u)K',H'=K'\cap\Ad(u\i)K=\Ad(u\i)H$.
Let $\Ps'_u=\Ps_u[[\d]]:\cd(V_{H,D})@>>>\cd(V_{H',D})$. We show:

(a) {\it If $A\in CS(V_{H,D})$ then $\Ps'_u(A)\in CS(V_{H',D})$.}
\nl
We have a commutative diagram in which the upper squares are cartesian and the left and
right vertical arrows are smooth with connected fibres:
$$\CD
V_{H,D}@<j<<\Xi_u@>h>>V_{H',D}\\
@Af'_3AA @Af'_4AA     @Af'_5AA   \\
V'_H@<j'<<\Xi''_u@>h'>>V'_{H'}\\
@Vf_3VV    @Vf_4VV    @Vf_5VV\\
V^u_H@<j^u<<\Xi'_u@>h^u>>V^u_{H'}   \endCD$$
Here

$V^u_H=\{(X,gU_X)\in V_{H,D};X=X_u\}$, $V^u_{H'}=\{(Y,gU_Y)\in V_{H',D};Y=Y_u\}$,

$\Xi'_u=\{(X,Y,g(U_X\cap U_Y))\in\Xi_u;X=X_u,Y=Y_u\}$,

$V'_H=G^0/U_{X_u}\T V^u_H$, $V'_{H'}=G^0/U_{Y_u}\T V^u_{H'}$,

$\Xi''_u=G^0/(U_{X_u}\cap U_{Y_u})\T\Xi'_u$,

$j^u,h^u$ are the restrictions of $j,h$,

$j'(x(U_{X_u}\cap U_{Y_u}),(X_u,Y_u,g(U_{X_u}\cap U_{Y_u}))=(xU_{X_u},(X_u,gU_{X_u}))$,

$h'(x(U_{X_u}\cap U_{Y_u}),(X_u,Y_u,g(U_{X_u}\cap U_{Y_u}))=(xU_{Y_u},(Y_u,gU_{Y_u}))$,

$f'_3(xU_{X_u},(X_u,gU_{X_u}))=(xX_ux\i,xgx\i U_{xX_ux\i}))$,

$f_3(xU_{X_u},(X_u,gU_{X_u}))=(X_u,gU_{X_u})$,

$f'_4(x(U_{X_u}\cap U_{Y_u}),(X_u,Y_u,g(U_{X_u}\cap U_{Y_u}))=
(xX_ux\i,xY_ux\i,xgx\i(U_{xX_ux\i}\cap U_{xY_ux\i}))$,

$f'_5(xU_{Y_u},(Y_u,gU_{Y_u}))=(xY_ux\i,xgx\i U_{xY_ux\i}))$,

$f_4(x(U_{X_u}\cap U_{Y_u}),(X_u,Y_u,g(U_{X_u}\cap U_{Y_u}))=
(X_u,Y_u,g(U_{X_u}\cap U_{Y_u}))$,

$f_5(xY_u,(Y_u,gU_{Y_u}))=(Y_u,gU_{Y_u})$.
\nl
We may identify $V^u_H=N_DX_u/U_{X_u}$, $V^u_{H'}=N_DY_u/U_{Y_u}$ in an obvious way. We can
find $C\in CS(N_DX_u/U_{X_u})$ such that $f'_3{}^\bst A=f_3^\bst C$. Hence 
$f'_3{}^*A[\dim X_u/U_{X_u}]=f_3^*C[\dim G/U_{X_u}]$ and $f'_3{}^*A=f_3^*A_1$ where 
$A_1=C[\a_u]$. Let $A_2=h^u_!j^{u*}A_1,A'_2=h_!j^*A$. Using 38.6(a) and the fact that $h^u$
is an isomorphism, we have
$$\align&f'_5{}^*A'_2=f'_5{}^*h_!j^*A=h'_!f'_4{}^*j^*A=h'_!j'{}^*f'_3{}^*A
=h'_!j'{}^*f_3^*A_1=h'_!f_4^*j^{u*}A_1\\&=h'_!f_4^*h^{u*}h^u_!j^{u*}A_1=
h'_!h'{}^*f_5^*h^u_!j^{u*}A_1=f_5^*h^u_!j^{u*}A_1[[-\d]]=f_5^*A_2[[-\d]].\endalign$$
Thus, $f_5^*A_2[[-\d]]=f'_5{}^*A'_2$. Hence 
$$f_5^\bst A_2[-\dim G^0/U_{Y_u}][[-\d]]=f'_5{}^\bst A'_2[-\dim Y_u/U_{Y_u}],
f_5^\bst A_2[-\a_u]=f'_5{}^\bst A'_2[[\d]].$$
Since $h^u,j^u$ are isomorphisms we see that $A_2[-\a_u]$ is a perverse sheaf. Hence so is
$f_5^\bst A_2[-\a_u]$. Hence $f'_5{}^\bst A'_2[[\d]]$ is perverse. Hence $A'_2[[\d]]$ is 
perverse and $h_!j^*A[[\d]]$ is perverse. To show that 
$h_!j^*A[[\d]]=A'_2[[\d]]\in CS(V_{H',D})$, it is enough to show that 
$A_2[-\a_u]\in CS(N_DY_u/U_{Y_u})$ or that $h^u_!j^{u*}C\in CS(N_DY_u/U_{Y_u})$.

Let $M$ be a common Levi subgroup of $X_u,Y_u$. Let $\tM=N_GX_u\cap N_GY_u\cap N_GM$. Then
$\tM$ is a reductive group with $\tM^0=M$ and $\tM^1=N_DX_u\cap N_DY_u\cap N_DM$ is a 
connected component of $\tM$. Moreover, the obvious maps $\tM^1@>>>N_DX_u/U_{X_u}$,
$\tM^1@>>>N_DY_u/U_{Y_u}$, $\tM^1@>>>(N_DX_u\cap N_DY_u)/(U_{X_u}\cap U_{Y_u})$ are 
isomorphisms (see 1.25). Hence the bottom row $V^u_H@<j^u<<\Xi'_u@>h^u>>V^u_{H'}$ of the
commutative diagram above  may be identified with $\tM^1@<<<\tM^1@>>>\tM^1$ where both maps
are the identity. Thus $h^u_!j^{u*}C\in CS(N_DY_u/U_{Y_u})$ follows immediately from 
$C\in CS(N_DX_u/U_{X_u})$. This proves (a).

We now show:

(b) {\it For $A$ as above we have $f_{H',J}\Ps'_u(A)=f_{H,J}(A)$.} 
\nl
Let $Q\in\cp_J$ be such that $X_u\sub Q\supset Y_u$. Let $M_1$ be the unique Levi of $Q$
such that $M\sub M_1$. Let $\tM_1=N_GQ\cap N_GM_1$. Then $\tM_1$ is a reductive group with
$\tM^0_1=M_1$ and $\tM^1_1=N_DQ\cap N_DM_1$ is a connected component of $\tM_1$. Moreover,
$X_u\cap\tM_1,Y_u\cap\tM_1$ are parabolic subgroups of $M_1$ with a common Levi, $M$. Let 
$C\in CS(N_DX_u/U_{X_u})=CS(\tM^1)$ be as above. We may assume that $A$ is simple so that 
$C$ is also simple. Using the proof of (a) and that of 38.2 we see that it is enough to 
verify that $\ind_{\tM^1}^{\tM^1_1}(C)$ (defined in terms of the parabolic $X_u\cap\tM_1$)
is isomorphic to $\ind_{\tM^1}^{\tM^1_1}(C)$ (defined in terms of the parabolic 
$Y_u\cap\tM_1$). Since $C$ is an admissible complex on $\tM^1$ (see 30.12), this follows 
from 27.2(d) which shows (for $G$ instead of $\tM_1$) that $\ind_{\tM^1}^{\tM^1_1}(C)$ can
be defined without reference to a choice of parabolic. This proves (b).

\proclaim{Proposition 38.8}Let $K,K',J$ be as in 38.4. Let $A\in CS(V_{K,D})$. We have
$$e_{K',J}f_{K,J}A\cong
\op_{u\in{}^K\WW^{K'}\cap\WW_J^\e}f_{K'\cap\Ad(u\i)K,K'}\Ps'_ue_{K\cap\Ad(u)K',K}A$$
in $CS(V_{K',D})$.
\endproclaim
We set $\fC=e_{K',J}f_{K,J}A\in\cd(V_{K,D})$. For any $u\in{}^K\WW^{K'}\cap\WW_J^\e$ we set
$\fC_u=f_{K'\cap\Ad(u\i)K,K'}\Ps'_ue_{K\cap\Ad(u)K',K}A\in\cd(V_{K',D})$. Assume that we 
can show that 
$$\fC\Bpq\{\fC_u;u\in{}^K\WW^{K'}\cap\WW_J^\e\}.$$
From the definition of $\Bpq$ (see 32.15) it would then follow that
$$\sum_i(-1)^i({}^pH^i(\fC))
=\sum_{u\in{}^K\WW^{K'}\cap\WW_J^\e}\sum_i(-1)^i({}^pH^i(\fC_u))$$
in the Grothendieck group of the category of perverse sheaves on $V_{K',D}$. From 
38.2(a), 38.3(a), 38.7(a) we see that $\fC,\fC_u\in CS(V_{K',D})$ (hence are perverse
sheaves); hence the previous equality implies that
$\fC=\sum_{u\in{}^K\WW^{K'}\cap\WW_J^\e}\fC_u$
in the Grothendieck group of the category of perverse sheaves on $V_{K',D}$. Since 
$\fC,\fC_u$ are semisimple perverse sheaves (being in $CS(V_{K',D})$) it follows that
$\fC\cong\op_{u\in{}^K\WW^{K'}\cap\WW_J^\e}\fC_u$, as desired.

Assume that $P\in\cp_K$ contains $X_u$ and $R\in\cp_{K'}$ contains $Y_u$ so that 
$\po(P,R)=u$. ($X_u,Y_u$ as in 38.6.) Let $Q\in\cp_J$ be such that $P\sub Q\supset R$. Let
$\b=\dim U_Q$. We have

$\fC=\ti{\fC}[\a+\a'-2\b](\a'-\b)$,

$\fC_u=\ti{\fC_u}[\a_u-\a](\a_u-\a)[\a_u-\a'][2\d](\d)$.
\nl
(Notation of 38.5.) Hence it is enough to show

$$\ti{\fC}[\a+\a'-2\b](\a'-\b)\Bpq\{\ti{\fC}_u[2\a_u-\a-\a'+2\d](\a_u-\a+\d);
u\in{}^K\WW^{K'}\cap\WW_J^\e\}$$
or that
$$\ti{\fC}\Bpq\{\ti{\fC}_u[2\a_u-2\a-2\a'+2\d-2\b](\a_u-\a-\a'+\d+\b);
u\in{}^K\WW^{K'}\cap\WW_J^\e\}.$$
By 38.5, it is enough to show that for any $u$ we have $\a_u-\a-\a'+\d+\b=-m_u$ or that
$$\dim U_P+\dim U_R-\dim(U_P\cap U_R)
=\dim U_{X_u}+\dim U_{Y_u}-\dim(U_{X_u}\cap U_{Y_u}).$$
It is enough to show that
$$\align&\dim\Lie U_P+\dim\Lie U_R-\dim(\Lie U_P\cap\Lie U_R)\\&=
\dim\Lie U_{X_u}+\dim\Lie U_{Y_u}-\dim(\Lie U_{X_u}\cap\Lie U_{Y_u})\endalign$$
or that $\dim(\Lie U_P+\Lie U_R)=\dim(\Lie U_{X_u}+\Lie U_{Y_u})$. We have 

$\Lie U_{X_u}=\Lie U_P+(\Lie P\cap\Lie U_R),\Lie U_{Y_u}=\Lie U_R+(\Lie R\cap\Lie U_P)$
\nl
hence $\Lie U_{X_u}+\Lie U_{Y_u}\sub\Lie U_P+\Lie U_R$. The opposite inclusion is clear 
since $\Lie U_P\sub\Lie U_{X_u},\Lie U_R\sub\Lie U_{Y_u}$. Thus we have 

$\Lie U_{X_u}+\Lie U_{Y_u}=\Lie U_P+\Lie U_R$. 
\nl
This completes the proof.

\subhead 38.9\endsubhead
For $J\sub\II$ such that $\e(J)=J$ let $\ck(V_{J,D})$ be the Grothendieck group of 
$CS(V_{J,D})$. Similarly let $\ck(D)$ be the Grothendieck group of $CS(D)$. For 
$A,A'\in CS(V_{J,D})$ we set $(A,A')=\dim\Hom(A,A')$. This induces a symmetric bilinear 
pairing $(,)_J:\ck(V_{J,D})\T\ck(V_{J,D})@>>>\ZZ$.

For $J\sub J'\sub\II$ such that $\e(J)=J,\e(J')=J'$, the functors
$f_{J,J'}:CS(V_{J,D})@>>>CS(V_{J',D})$ and $e_{J,J'}:CS(V_{J',D})@>>>CS(V_{J,D})$ are 
compatible with direct sums hence they induce homomorphisms
$\ck(V_{J,D})@>>>\ck(V_{J',D})$, $\ck(V_{J',D})@>>>\ck(V_{J,D})$ denoted again by
$f_{J,J'},e_{J,J'}$. From 30.5 we see that 

(a) $(e_{J,J'}A',A)_J=(A',f_{J,J'}A)_{J'}$
\nl
for $A\in\ck(V_{J,D}),A'\in\ck(V_{J,D})$.

In the setup of 38.7, for $u\in{}^K\WW^{K'}\in\WW_J^\e$, the functor 

$\Ps'_u:CS(V_{K\cap\Ad(u)K',D})@>>>CS(V_{K'\cap\Ad(u\i)K,D}$
\nl
is compatible with direct sums hence induces a homomorphism 

$\ck(V_{K\cap\Ad(u)K',D})@>>>\ck(V_{K'\cap\Ad(u\i)K,D}$ 
\nl
denoted again by $\Ps'_u$. Below we shall need the following identity:
$$\sum\Sb K';K'\sub J\\\e(K')=K'\eSb\sum\Sb u\in{}^K\WW^{K'}\cap\WW_J^\e\\K\cap\Ad(u)K'=H
\eSb(-1)^{|K'_\e|}=(-1)^{|H_\e|}\tag b$$
for any $H\sub K\sub J\sub\II$ such that $\e(H)=h,\e(K)=K,\e(J)=J$. In the case where 
$J=\II,\e=1$ this is proved in \cite{\CU, 2.5}; the general case can be reduced to this 
special case by replacing $\WW$ by $\WW_J^\e$ which is itself a Weyl group with simple 
reflections in bijection with $J_\e$.

\subhead 38.10\endsubhead
For $J\sub\II$ such that $\e(J)=J$ we define a homomorphism
$\dd_J:\ck(V_{J,D})@>>>\ck(V_{J,D})$ by
$$\dd_JA=\sum_{K;K\sub J,\e(K)=K}(-1)^{|K_\e|}f_{K,J}e_{K,J}A.\tag a$$
Using 38.9(a) we see that for $A,A'\in\ck(V_{J,D})$ we have
$$(\dd_J(A),A')_J=(A,\dd_J(A')_J.\tag b$$
We show that, for $K\sub J\sub\II$ such that $\e(K)=K,\e(J)=J$ and $A\in\ck(V_{K,D})$, we
have
$$\dd_Jf_{K,J}A=f_{K,J}\dd_KA.\tag c$$
Using 38.8, 38.7(b), 38.9(b), 38.1(a) we have
$$\align&\dd_Jf_{K,J}A=\sum_{K';K'\sub J,\e(K')=K'}(-1)^{|K'_\e|}f_{K',J}e_{K',J}f_{K,J}A
\\&=\sum\Sb K';K'\sub J\\\e(K')=K'\\u\in{}^K\WW^{K'}\cap\WW_J^\e\eSb
(-1)^{|K'_\e|}f_{K',J}f_{K'\cap\Ad(u\i)K,K'}\Ps'_ue_{K\cap\Ad(u)K',K}A\\&
=\sum\Sb H;H\sub J,\e(H)=H\\K'\\ 
u\in{}^K\WW^{K'}\cap\WW_J^\e\\K'\sub J,\e(K')=K'\\
K\cap\Ad(u)K'=H\eSb(-1)^{|K'_\e|}f_{\Ad(u\i)H,J}\Ps'_ue_{H,K}A\\&
=\sum\Sb H;H\sub J,\e(H)=H\\K';u\in{}^K\WW^{K'}\cap\WW_J^\e\\K'\sub J,\e(K')=K'\\
K\cap\Ad(u)K'=H\eSb(-1)^{|K'_\e|}f_{H,J}e_{H,K}A\\&
=\sum\Sb H;H\sub J,\e(H)=H\\K';K'\sub J,\e(K')=K'\\
u\in{}^K\WW^{K'}\cap\WW_J^\e\\K\cap\Ad(u)K'=H\eSb(-1)^{|K'_\e|}f_{H,J}e_{H,K}A
=\sum_{H;H\sub K,\e(H)=H}(-1)^{|H_\e|}f_{H,J}e_{H,K}A\\&
=\sum_{H;H\sub K,\e(H)=H}(-1)^{|H_\e|}f_{K,J}f_{H,K}e_{H,K}A=f_{K,J}\dd_KA\endalign$$
and (c) is proved.

We show that, for $J\sub\II$ such that $\e(J)=J$ and $A\in\ck(V_{J,D})$, we have
$$\dd_J\dd_JA=A.\tag d$$
Using (c), 38.1(a), we have
$$\align&\dd_J\dd_JA=\sum_{K;K\sub J,\e(K)=K}(-1)^{|K_\e|}\dd_Jf_{K,J}e_{K,J}A\\&=
\sum_{K;K\sub J,\e(K)=K}(-1)^{|K_\e|}f_{K,J}\dd_Ke_{K,J}A\\&
=\sum_{K,K';K'\sub K\sub J,\e(K)=K,\e(K')=K'}(-1)^{|K_\e|}(-1)^{|K'_\e|}
f_{K,J}f_{K',K}e_{K',K}e_{K,J}A\\&=\sum_{K,K';K'\sub K\sub J,\e(K)=K,\e(K')=K'}
(-1)^{|K_\e|}(-1)^{|K'_\e|}f_{K',J}e_{K',J}A\\&=\sum_{K';K'\sub J,\e(K')=K'}
(-1)^{|K'_\e|}\sum_{K;K'\sub K\sub J,\e(K)=K}(-1)^{|K_\e|}f_{K',J}e_{K',J}A\\&
=\sum_{K';K'\sub J,\e(K')=K'}(-1)^{|K'_\e|}\d_{K',J}(-1)^{K'_\e}f_{K',J}e_{K',J}A
=f_{J,J}e_{J,J}A=A\endalign$$
and (d) is proved.

We show that, for $J\sub\II$ such that $\e(J)=J$ and $A,A'\in\ck(V_{J,D})$, we have
$$(\dd_JA,\dd_JA')_J=(A,A')_J.\tag e$$
Using (b),(d) we have $(\dd_JA,\dd_JA')_J=(A,\dd_J\dd_JA')_J=(A,A')_J$ as desired.

\subhead 38.11\endsubhead
We write $\dd$ instead of $\dd_{\II}$ and $(,)$ instead of $(,)_\II$. We call $\dd$ the
duality operator on character sheaves. If $A$ is a character sheaf on $D=V_{\II,D}$, then
$(A,A)=1$ (where $A$ is regarded as an element of $\ck(D)$) hence, by 38.10(e), we have 
$(\dd(A),\dd(A))=1$. Since $\dd(A)$ is a $\ZZ$-linear combination of isomorphism classes of
character sheaves (which form an orthonormal basis of $\ck(D)$ for $(,)$) it follows that 
$\dd(A)=\pm A'$ where $A'$ is a well defined character sheaf on $D$ (up to isomorphism). 
The sign can be described as follows. By 30.6(d) we can find a parabolic
$P_0$ of $G^0$ such that $N_DP_0\ne\em$ and a cuspidal character sheaf $A_0$ on 
$D_0:=N_DP_0/U_{P_0}$ such that $A$ is a direct summand of $\ind_{D_0}^D(A_0)$. We have 
$P_0\in\cp_J$ where $J\sub\II,\e(J)=J$. Then
$$\dd(A)=(-1)^{|J_\e|}A'.\tag a$$
Indeed let $A_0^\flat\in CS(V_{D,J})$ be the perverse sheaf corresponding to $A_0$ as in
30.3. Then $e_{J',J}A_0^\flat=0$ for any $J'\subsetneq J$ such that $\e(J')=J'$, see 38.2.
Hence $\dd_J(A_0^\flat)=(-1)^{|J_\e|}A_0^\flat$. Now $A$ is a direct summand of 
$f_{J,\II}A_0^\flat$ and by 38.10(c) we have 
$$\dd f_{J,\II}A_0^\flat=f_{J,\II}\dd_JA_0^\flat=(-1)^{|J_\e|}f_{J,\II}A_0^\flat.$$
In $\ck(D)$ we have $f_{J,\II}A_0^\flat=\sum_{k=1}^{k_0}n_kA_k$ where $A_k$ are distinct
character sheaves on $D$, $n_k\in\ZZ_{>0}$ and $A_1=A$. We have $\dd(A_k)=\io_kA'_k$ 
where $\io_k=\pm1$ and $A'_k$ are distinct character sheaves on $D$, $n_k\in\ZZ_{>0}$ and
$A_1=A$. We see that $\sum_kn_k\io_kA'_k=(-1)^{|J_\e|}\sum_kn_kA_k$. Since $\{A_k\}$ and 
$\{A'_k\}$ are parts of the same basis of $\ck(D)$ we see that $\io_k/(-1)^{|J_\e|}>0$ 
for any $k$. Hence $\io_k=(-1)^{|J_\e|}$ for any $k$. In particular this holds for $k=1$
and (a) follows.

Note that, by 38.10(d), $A\m A'$ is an involution of the set of isomorphism classes of 
character sheaves on $D$. 

\subhead 38.12\endsubhead
The definition of the duality operator for character sheaves in 38.11 is entirely similar
to that of a duality operator for representations of a reductive group over a finite field
(given again by an alternating sum of compositions of a parabolic restriction and a
parabolic induction) which was found by the author in 1977, who conjectured that it takes 
irreducibles to $\pm$ irreducibles and is involutive. In 1977 I communicated this
conjecture to C.W.Curtis and N.Kawanaka, see \cite{\KA, p.412}; the conjecture was proved 
in \cite{\CU},\cite{\AL} and in \cite{\KA}.

\head 39. Quasi-rationality\endhead
\subhead 39.1\endsubhead
The main result of this section is Proposition 39.7 which gives a quasi-rationality
property of representations of certain extensions of a Weyl group. This is needed to prove
a key property of character sheaves (Corollary 39.8).

Let $W,I$ be a Weyl group ($I$ is the set of simple reflections). We have canonically 
$W=\prod_{j\in J}W_j$ where $W_j$ is an irreducible Weyl group (with set of simple 
reflections $I_j=I\cap W_j$). We identify $W_j$ with a subgroup of $W$. Let 
$\Rf(W)=\cup_{w\in W}wIw\i$. We have $\Rf(W)=\sqc_{j\in J}\Rf(W_j)$. For any $j\in J$ the 
set $\Rf(W_j)$ is a single $W_j$-conjugacy class if $W_j$ is of type $A,D$ or $E$; it is a
union of two $W_j$-conjugacy classes, otherwise. A subset $\cx$ of $\Rf(W)$ is said to be 
{\it special} if $\cx=\sqc_{j\in J}\cx_j$ where $\cx_j$ is a $W_j$-conjugacy class in 
$\Rf(W_j)$. 
Clearly, a special subset of $\Rf(W)$ exists; we assume that a special subset $\cx$ of 
$\Rf(W)$ is given. Let $\Aut(W,I,\cx)$ be the group of automorphisms of $W$ that preserve 
$I$ and $\cx$.

\proclaim{Lemma 39.2} Let $\G$ be a finite group. Let $\g\m\r_\g$ be a homomorphism 
$\G@>>>\Aut(W,I,\cx)$. Let $E$ be a simple $\QQ[W]$-module such that 
$\tr(\r_\g(w),E)=\tr(w,E)$ for any $\g\in\G,w\in W$. Assume that either

(i) $|J|=1$, or 

(ii) $\G$ is an extension of a cyclic group by a cyclic group.
\nl
Then there exists a homomorphism $\G@>>>\Aut_\QQ(E)$, $\g\m t_\g$ such that 
$t_\g(w(e))=\r_\g(w)t_\g(e)$ for any $\g\in\G,w\in W,e\in E$.
\endproclaim
$\G$ acts on $J$ by $\g:j\m\g(j)$ where $W_{\g(j)}=\r_\g(W_j)$. We may identify
$E=\ot_{j\in J}E_j$ as $\QQ[W]$-modules where $E_j$ is a simple $\QQ[W_j]$-module for any
$j\in J$. From our assumption we see that for any $j\in J,\g\in\G$ there exists a
$\QQ$-linear isomorphism 

(a) $H_j^\g:E_j@>\si>>E_{\g(j)}$ with $H_j^\g(w_je)=\r_\g(w_j)H_j^\g(e)$ for all 
$e\in E_j,w_j\in W_j$;
\nl
moreover, $H_j^\g$ is unique up to multiplication by an element of $\QQ^*$. It follows that

(b) $H_{\g(j)}^{\g'}H_j^\g\in\QQ^*H_j^{\g'\g}$ for $j\in J$ and $\g,\g'\in\G$.
\nl
For $\g\in\G$ we define $\tit_\g\in\Aut_\QQ(E)$ by

(c) $\tit_\g(\ot_je_j)=\ot_je'_j$ where $e_j\in E_j$ and 
$e'_j=H_{\g\i(j)}^\g(e_{\g\i(j)})\in E_j$.
\nl
From definitions we have $\tit_\g(w(e))=\r_\g(w)\tit_\g(e)$ for any $w\in W,e\in E$;
moreover for $\g,\g'\in\G$, the maps $\tit_{\g'\g},\tit_{\g'}\tit_\g$ are equal up to a 
factor in $\QQ^*$. Thus the maps $\tit_\g$ provide a homomorphism $\G@>>>PGL(E)$ rather 
than a homomorphism $\G@>>>GL(E)$.

We prove the lemma in the setup of (i). Replacing $\G$ by its image under $\g\m\r_g$ we may
assume that $\G\sub\Aut(W,I,\cx)$ and $\g\m\r_\g$ is the inclusion. We form the semidirect
product $W\G$ with $W$ normal. It is enough to show that $E$ 
extends to a $W\G$-module. We may assume that $\G\ne\{1\}$. Then $W$ is of type 
$A_n(n\ge 2)$, $D_n$ or $E_6$. If $\G=\Aut(W,I,\cx)$ then $W\G$ is itself a Weyl group, of
type $A_n\T A_1,D_{2n+1}\T A_1,B_{2n},F_4,E_6\T A_1$ for $W$ of type 
$A_n(n\ge 2),D_{2n+1}(n\ge 2),D_{2n}(n\ge 3),D_4,E_6$ respectively and the desired result
follows easily from the known properties of representations of such Weyl groups (in 
particular, from their rationality); the same applies if $W$ is of type $D_4$ and 
$|\G|=2$ (in this case, $W\G$ is a Weyl group of type $B_4$). In the only remaining case 
($W$ of type $D_4,|\G|=3$), the result follows by an argument in \cite{\CRG, 3.2}.

Next we prove the lemma in the setup of (ii). Now $\G$ has two generators $a,c$ and 
relations 

$a^M=1$, $c^N=a^u$, $cac\i=a^k$
\nl
where $M,N,k$ are integers $\ge 1$  such that $k^N=1\mod M$ and $u\in\NN$ satisfies
$uk=k\mod M$.

We consider separately 3 cases in increasing order of generality.

{\it Case 1.} Assume that $J$ is a single orbit of $a:J@>>>J$. We may identify $J=\ZZ/m$
so that $a(i)=i+1,c(i)=ki-r$ for $i\in\ZZ/m$; here $r\in\NN$ is independent of $i$. Since
$a^M=1:J@>>>J$ we see that $m$ divides $M$. For any $i$ we have 
$i+u=c^N(i)=k^Ni-(1+k+k^2+\do+k^{N-1})r$ in $\ZZ/m$. Since $k^N=1\mod M$ (hence 
$k^N=1\mod m$) we have $u+(1+k+k^2+\do+k^{N-1})r=fm$ for some $f\in\NN$. We set $c'=a^rc$.
Then $c'(i)=ki$ for $i\in\ZZ/m$. Let $\G'$ be the subgroup of $\G$ generated by $a^m,c'$. 
In $\G'$ we have 

$(a^m)^{M/m}=1$, $c'{}^N=(a^m)^f$, $c'a^mc'{}\i=(a^m)^k$.
\nl
Since $a^m(0)=0,c'(0)=0$, the action of $\G$ on $W$ restricts to an action of $\G'$ on 
$W_0$. Using the lemma (setup of (i)) for $W_0,E_0,\G'$ instead of $W,E,\G$ we obtain a
homomorphism $\G'@>>>\Aut_\QQ(E_0)$ such that, denoting by $s_0,s_1$ the images of $a^m,c'$
under this homomorphism, we have $s_0^{M/m}=1$, $s_1^N=s_0^f$, $s_1s_0s_1\i=s_0^k$, 
$s_0(w(e))=\r_a^m(w)s_0(e)$, $s_1(w(e))=\r_{c'}(w)s_1(e)$ for any $w\in W_0,e\in E_0$. 

For $i\in\ZZ/m$, let $H_i^a:E_i@>\si>>E_{i+1}$ be as in (a). Then 
$s_2:=H_{m-1}^a\do H_1^aH_0^a:E_0@>\si>>E_0$ satisfies $s_2(w(e))=\r_a^m(w)s_2(e)$ for any
$w\in W_0,e\in E_0$. Hence $s_2\i s_0(w(e))=w(s_2\i s_0(e))$ for any $w\in W_0$, 
$e\in E_0$. By the absolute irreducibility of the $W_0$-module $E_0$ we see that 
$s_2\i s_0:E_0@>>>E_0$ is a $\QQ^*$-multiple of the identity map. Hence, replacing $H_0^a$
by a $\QQ^*$-multiple, we may assume that $s_2=s_0$ that is, $H_{m-1}^a\do H_1^aH_0^a=s_0$.
For any $h\in\NN$ we set

$b_h=H_{hk-1}^aH_{hk-2}^a\do H_1^aH_0^as_1(H_{h-1}^a\do H_1^aH_0^a)\i:E_h@>\si>>E_{kh}$
\nl
(there are $hk$ factors to the left of $s_1$ and $h$ factors to the right of $s_1$). For
$h\ge 1$ we have 

(d) $b_h=H_{hk-1}^aH_{hk-2}^a\do H_{hk-k}^ab_{h-1}(H_{h-1}^a)\i$.
\nl
We show that $b_{h+m}=b_h$ for $h\in\NN$. We argue by induction on $h$. Assume first that 
$h=0$. We must show that 
$H_{mk-1}^aH_{mk-2}^a\do H_1^aH_0^as_1(H_{m-1}^a\do H_1^aH_0^a)\i=s_1$; this is a
reformulation of the already known equality $s_0^ks_1s_0\i=s_1$. Assume next that $h\ge 1$.
Using (d) twice and the induction hypothesis we have
$$\align&b_{h+m}=H_{hk+mk-1}^aH_{hk+mk-2}^a\do H_{hk+mk-k}^ab_{h+m-1}(H_{h+m-1}^a)\i\\&=
H_{hk-1}^aH_{hk-2}^a\do H_{hk-k}^ab_{h-1}(H_{h-1}^a)\i=b_h,\endalign$$
as desired. We see that $b_h$ depends only on the image of $h$ in $\ZZ/m$.

We set $k=k^{N-1}$. Then $kk'=1\mod M$. Define $t_a,t_{c'}\in\Aut_\QQ(E)$ by 
$t_a(\ot_{i\in\ZZ/m}e_i)=\ot_{i\in\ZZ/m}e'_i$ where $e_i\in E_i$ and 
$e'_i=H_{i-1}^a(e_{i-1})\in E_i$,

$t_{c'}(\ot_{i\in\ZZ/m}e_i)=\ot_{i\in\ZZ/m}e''_i$ where $e_i\in E_i$ and 
$e''_i=b_{ik'}(e_{ik'})\in E_{ik'k}=E_i$.
\nl
From definitions we have $t_a(w(e))=\r_a(w)t_a(e),t_{c'}(w(e))=\r_{c'}(w)t_{c'}(e)$ for
any $w\in W,e\in E$. 

We have $t_{c'}t_at_{c'}\i=t_a^k$. This follows from the identity

$b_{hk'}H_{hk'-1}^ab_{hk'-1}\i=H_{h-1}^aH_{h-2}^a\do H_{h-k}^a:E_{h-k}@>>>E_h$
\nl
for $h=0,1,\do,m-1$ (here $b_{-1}$ is taken to be $b_{m-1}$); an equivalent identity is

$b_iH_{i-1}^ab_{i-1}\i=H_{ik-1}^aH_{ik-2}^a\do H_{ik-k}^a:E_{ik-k}@>>>E_{ik}$
\nl
for $i=0,1,\do,m-1$, which is the same as (d).

We have $t_a^M=1$. This follows from the identity
$H_{i+M-1}^a\do H_{i+1}^aH_i^a=1:E_i@>>>E_i$ for $i=0,1,\do,m-1$ which is equivalent to the
known equality $s_0^{M/m}=1$.

We show that 

(e) $t_{c'}^N=t_a^{mf}$. 
\nl
This follows from the identity

$b_{ik'}b_{ik'{}^2}\do b_{ik'{}^N}=H_{i-1}^aH_{i-2}^a\do H_{i-mf}^a$ for $i=0,1,\do,m-1$ 
\nl
or equivalently

$b_{ik^{N-1}}b_{ik^{N-2}}\do b_i=H_{i-1}^aH_{i-2}^a\do H_{i-mf}^a$ for $i=0,1,\do,m-1$.
\nl
This is the same as

$H_{ik^N-1}^aH_{ik^N-2}^a\do H_1^aH_0^as_1^N(H_{i-1}^a\do H_1^aH_0^a)\i=
H_{i-1}^aH_{i-2}^a\do H_{i-mf}^a$
\nl
(with $ik^N$ factors to the left of $s_1^N$). We have $ik^N-i=iMl$ for some $l\in\NN$. From
$s_0^{M/m}=1$ we deduce $s_0^{iMl/m}=1$ hence $H_{ik^N-1}^aH_{ik^N-2}^a\do H_i^a=1$ (with 
$ik^N-i$ factors). Since $s_1^N=s_0^f$, it remains to show  

$H_{i-1}^a\do H_1^aH_0^as_0^f(H_{i-1}^a\do H_1^aH_0^a)\i=H_{i-1}^aH_{i-2}^a\do H_{i-mf}^a$.
\nl
If $f=0$ this is obvious. Assume now that $f\ge 1$. Then $i-mf<0$ and we see that it is
enough to show

$s_0^f(H_{i-1}^a\do H_1^aH_0^a)\i=H_{-1}^aH_{-2}^a\do H_{i-mf}^a$
\nl
that is,

$s_0^f=H_{-1}^aH_{-2}^a\do H_{i-mf}^a(H_{i-1}^a\do H_1^aH_0^a)$.
\nl
This can be rewritten in the form

$s_0^f=H_{-1}^aH_{-2}^a\do H_{i-mf}^aH_{i-mf-1}^a\do H_{1-mf}^aH_{-mf}^a$
\nl
which follows from $s_0=H_{m-1}^a\do H_1^aH_0^a$. This proves (e).

We set $t_c=t_a^{-r}t_{c'}\in\Aut_\QQ(E)$. From definitions we have 
$t_c(w(e))=\r_c(w)t_c(e)$ for any $w\in W,e\in E$. 

We have $t_ct_at_c\i=t_a^k$. This follows from the identity $t_{c'}t_at_{c'}\i=t_a^k$.

We have $t_c^N=t_a^u$. Indeed, 
$$(t_a^{-r}t_{c'})^N=t_a^{-r(1+k+k^2+\do+k^{N-1})}t_{c'}^N=t_a^{-mf+u}t_a^{mf}=t_a^u.$$
We see that $t_a,t_c\in\Aut_\QQ(E)$ satisfy the relations of $\G$ hence they define a
homomorphism $\G@>>>\Aut_\QQ(E)$. This has the required properties.

{\it Case 2.} Assume that $J$ is a single $\G$-orbit. Let $\la a\ra$ be the subgroup of 
$\G$ generated by $a$. If $X$ is an $\la a\ra$-orbit in $J$ then cX is again an 
$\la a\ra$-orbit. (We must show that, if $j\in X$ and $i\ge 1$, then $ca^ij,cj$ are in the
same $\la a\ra$-orbit. But $ca^ij=a^{ik}cj$.) Hence $c^hX$ is an $\la a\ra$-orbit in $J$ 
for any $h\in\NN$. We can find an integer $z\ge 1$ such that $X,cX,...,c^{z-1}X$ are 
distinct and $c^zX=X$. (Clearly, $z$ is a divisor of $N$.) Hence the notation $X_h=c^hX$ 
for $h\in\ZZ/z$ is meaningful. Now $X_0,X_1,\do,X_{z-1}$ are distinct. The union 
$\cup_{h\in[0,z-1]}X_h$ is $c$-stable and $a$-stable hence it is equal to $J$ (by our 
assumption on $J$). We see that $|J|=|X|z$. For $h\in\ZZ/z$ we set 
$W^h=\prod_{j\in X_h}W_j$, $E^h=\ot_{j\in X_h}E_j$. We have naturally 
$W=\prod_{h\in[0,z-1]}W^h$, $E=\ot_{h\in[0,z-1]}E^h$. Let $\G''$ be the subgroup of $\G$
generated by $a,c^z$. In $\G''$ we have 

$a^M=1$, $(c^z)^{N/z}=a^u$, $c^zac^{-z}=a^{k^z}$.
\nl
Since $a(X_0)=X_0,c^z(X_0)=X_0$, the action of $\G$ on $W$ restricts to an action of $\G''$
on $W^0$. Using case 1 for $W^0,E^0,\G''$ instead of $W,E,\G$ we obtain a homomorphism 
$\G''@>>>\Aut_\QQ(E^0)$ such that, denoting by $S_0,S_1$ the images of $a,c^z$ under this 
homomorphism, we have $S_0^M=1$, $S_1^{N/z}=S_0^u$, $S_1S_0S_1\i=S_0^{k^z}$, 
$S_0(w(e))=\r_a(w)S_0(e)$, $S_1(w(e))=\r_c^z(w)S_1(e)$ for any $w\in W^0,e\in E^0$. 

For $h\in\ZZ/z$ there exists a $\QQ$-linear isomorphism $K_h:E^h@>\si>>E^{h+1}$ with 
$K_h(we)=\r_c(w)K_h(e)$ for all $e\in E^h,w\in W^h$ (note that $\r_c(w)\in W^{h+1}$);
moreover, $K_h$ is unique up to multiplication by an element of $\QQ^*$. (For example we
can take $K_h$ of the form $\ot_{j\in X_h}H_j^c$ where $H_j^c:E_j@>\si>>E_{cj}$ are as in
(a).) Then $S_2:=K_{z-1}\do K_1K_0:E^0@>\si>>E^0$ satisfies $S_2(w(e))=\r_c^z(w)S_2(e)$ for
any $w\in W^0,e\in E^0$. Hence $S_2\i S_1(w(e))=w(S_2\i S_1(e))$ for any 
$w\in W^0,e\in E^0$. By the absolute irreducibility of the $W^0$-module $E^0$ we see that
$S_2\i S_1:E^0@>>>E^0$ is a $\QQ^*$-multiple of the identity map. Hence, replacing $K_0$ by
a $\QQ^*$-multiple, we may assume that $S_2=S_1$ that is, $K_{z-1}\do K_1K_0=S_1$.

For any $h\in\NN$ we define $\b_{-h}\in\Aut_\QQ(E^{-h})$ by 

$\b_{-h}=K_{-h}\i\do K_{-2}\i K_{-1}\i S_0^{k^h}K_{-1}K_{-2}\do K_{-h}$
\nl
(for $h=0$ this is interpreted as $\b_0=S_0$). We show that $\b_{-h-z}=\b_{-h}$. We argue 
by induction on $h$. For $h=0$ we must verify that 

$K_{-z}\i\do K_{-2}\i K_{-1}\i S_0^{k^z}K_{-1}K_{-2}\do K_{-z}=S_0$;
\nl
this follows from the known equality $S_1\i S_0^{k^z}S_1=S_0$. For $h\ge 1$ we have using 
the induction hypothesis

$\b_{-h-z}=K_{-h-z}\i\b_{-h-z+1}^kK_{-h-z}=K_{-h}\i\b_{-h+1}^kK_{-h}=\b_{-z}$,
\nl
We see that $\b_{-h}$ depends only on the image of $h$ in $\ZZ/z$. Define 
$t_a,t_c\in\Aut_\QQ(E)$ by 

$\t_a(\ot_{h\in\ZZ/z}e_h)=\ot_{h\in\ZZ/z}e'_h$ where $e_h\in E^h$ and 
$e'_h=\b_h(e_h)\in E^h$,

$\t_c(\ot_{h\in\ZZ/z}e_h)=\ot_{h\in\ZZ/z}e''_h$ where $e_h\in E^h$ and 
$e''_h=K_{h-1}(e_{h-1})\in E^h$.
\nl
From definitions we have $t_a(w(e))=\r_a(w)t_a(e),t_c(w(e))=\r_c(w)t_c(e)$ for any 
$w\in W,e\in E$. 

We have $t_ct_at_c\i=t_a^k$. This follows from the identity
$\b_{h-1}=K_{h-1}\i\b_h^kK_{h-1}$ for $h\in\ZZ/z$.

We have $t_a^M=1$. This follows from the identity $\b_{-h}^M=1$ for any $h\in\NN$. An
equivalent statement is $K_{-h}\i\do K_{-2}\i K_{-1}\i S_0^{k^hM}K_{-1}K_{-2}\do K_{-h}=1$
which follows from $S_0^M=1$.

We show that $t_c^N=t_a^u$. It is enough to show that
$K_{-h-1}K_{-h-2}\do K_{-h-N}=b_{-h}^u:E_{-h}@>>>E_{-h}$ for $h=0,1,\do,z-1$. An equivalent
statement is
$$(K_{-h}\i\do K_{-1}\i)(K_{-1}K_{-2}\do K_{-N})(K_{-1}\do K_{-h})=
K_{-h}\i\do K_{-1}\i S_0^{k^hu}K_{-1}\do K_{-h}.$$
The left hand side is $(K_{-h}\i\do K_{-1}\i)S_1^{N/z}(K_{-1}\do K_{-h})$. It is enough to
show that $S_1^{N/z}=S_0^{k^hu}$ or that $S_0^u=S_0^{k^hu}$ for $h=0,1,\do,z-1$. Since
$S_0^M=1$ it is enough to show that $u=k^hu\mod M$. This follows from $uk=u\mod M$.

We see that $t_a,t_c\in\Aut_\QQ(E)$ satisfy the relations of $\G$ hence they define a
homomorphism $\G@>>>\Aut_\QQ(E)$. This has the required properties.

{\it Case 3.} We now consider the general case. For any $\G$-orbit $Y$ on $J$ we set
$W^Y=\prod_{j\in Y}W_j$, $E^Y=\ot_{j\in Y}E_j$. We have naturally $W=\prod_YW^Y$, 
$E=\ot_YE^Y$ where $Y$ runs over the $\G$-orbits in $J$. Now the $\G$-action on $W$
restricts to a $\G$-action on $W^Y$ for each $Y$. Using case 2 for $W^Y,E^Y,\G$ instead of
$W,E,\G$ we obtain for each $Y$ a homomorphism $\G@>>>\Aut_\QQ(E^Y)$. We define a
homomorphism $\G@>>>\Aut_\QQ(E)$ by $\g:\ot_Ye_Y\m\ot_Y(\g(e_Y))$; here $e_Y\in E^Y$. This
has the required properties.

This completes the proof in the setup of (ii). The lemma is proved.

\subhead 39.3\endsubhead
Let $\G$ be a finite group. Assume that $\G$ is a semidirect product of a normal subgroup
$\G'$ with a cyclic group of order $n$ with generator $b$. Let $\bfU$ be an algebraic 
closed field of characteristic $0$. Let $E$ be a simple $\bfU[\G]$-module which is 
isotypical as a $\bfU[\G']$-module. We show that

(a) {\it $E$ is simple as a $\bfU[\G']$-module.}
\nl
There exists a simple $\bfU[\G']$-module $E'$ such that, setting 
$V=\Hom_{\bfU[\G']}(E',E)$ we have $V\ot E'@>\si>>E,f\ot e'\m f(e')$. For any $\g'\in\G'$ 
we have $\tr(b\g'b\i,E')=\tr(\g',E')$. Hence there exists $\x\in\Aut_{\bfU}(E')$ such that
$\x(\g'e')=(b\g'b\i)\x(e')$ for $e'\in E',w'\in\G'$. Then $\x^n\in\Aut_{\bfU}(E')$ commutes
with the $\G'$-action hence it is a scalar on $E'$. Replacing $\x$ by an $\bfU^*$-multiple
we can assume that $\x^n=1$. Then $E'$ becomes a $\bfU[\G]$-module with $b$ acting as $\x$.
Define $\et:V@>>>V$ by $f\m\et(f)$ where $\et(f)(e')=bf(\x\i(e'))$ for $e'\in E'$. We 
regard $V$ as a $\bfU[\G]$-module in which $\G'$ acts trivially and $b$ acts as $\et$. Then
the isomorphism $V\ot E'@>\si>>E$ considered above is an isomorphism of $\bfU[\G]$-modules.
Since $E$ is simple it follows that $V$ is a simple $\bfU[\G]$-module. Since $\G$ acts on 
$V$ through a cyclic quotient, we see that $\dim V=1$. It follows that $E\cong E'$ as a 
$\bfU[\G']$-module; (a) follows. 

\subhead 39.4\endsubhead
Let $W,I,\cx$ be as in 39.1. Let $\G$ be a finite group with generators $a,c$ and relations
$a^M=1$, $c^N=1$, $cac\i=a^k$ where $M,N,k$ are integers $\ge 1$ such that $k^N=1\mod M$. 
Let $\la a\ra$ (resp. $\la c\ra$) be the subgroup of $\G$ generated by $a$ (resp. $c$). Let
$\g\m\r_\g$ be a homomorphism $\G@>>>\Aut(W,I,\cx)$. We form the semidirect product 
$W\G,W\la a\ra$ with $W$ normal. Note that $W\la a\ra$ is a subgroup of $W\G$. 

\proclaim{Lemma 39.5}In the setup of 39.4 let $E$ be a simple $\bfU[W\G]$-module. Assume
that either $cxc\i=x$ for any $x\in\la a\ra$, or $cxc\i=x\i$ for any $x\in\la a\ra$. Let 
$x\in\la a\ra$. There exists $\z$, a root of $1$ in $\bfU$, such that $\tr(wxc,E)\in\z\ZZ$
for any $w\in W$.
\endproclaim
Let $W'=W\la a\ra$. We can write canonically $E=\op_{t\in T}E^t$ where $E^t$ are isotypical
$\bfU[W']$-modules. Now $\G$ acts on $T$ by $\g E^t=E_{\g(t)}$. This action is transitive
since $E$ is simple as a $\bfU[W\G]$-module. If $c(t)\ne t$ for any $t\in T$ then for any
$w\in W$, $wxc:E@>>>E$ permutes the summands $E^t$ and no summand is stable hence 
$\tr(wxc,E)=0$. In this case the lemma is clear. Thus we may assume that $c(t)=t$ for some
$t\in T$. Then $cE^t=E^t$ hence $E^t$ is a $\bfU[W\G]$-submodule of $E$. Since $E$ is 
simple, we have $E=E^t$. Thus, $E$ is isotypical as a $\bfU[W']$-module. Using 39.3(a) for
$W\G,W',\la c\ra$ instead of $\G,\G',C$ we see that $E$ is simple as a $\bfU[W']$-module. 

We can write canonically $E=\op_{h\in H}E_h$ where $E_h$ are isotypical 
$\bfU[W]$-modules. Now $\G$ acts on $H$ by $\g E_h=E_{\g(h)}$. This action is transitive
since $E$ is simple as a $\fU[W\G]$-module. The restriction of this action to $\la a\ra$ is
also transitive since $E$ is simple as a $\bfU[W']$-module. For $h\in H$ let $\la a\ra_h$ 
be the stabilizer of $h$ in $\la a\ra$. Let $W'_h=W\la a\ra_h\sub W'$. Then $E_h$ is a 
$W'_h$-submodule of $E$ and the $W'$-module $E$ is induced by the $W'_h$-module $E_h$. 
Since $E$ is simple as a $W'$-module it follows that $E_h$ is simple as a $W'_h$-module. 
Using 39.3(a) for $W'_h,W,E_h$ instead of $\G,\G',E$, we see that $E_h$ is simple as a 
$W$-module.

For $h\in H$ let $\G_h=\{\g\in\G;\g(h)=h\}$. Then $E_h$ is a $W\G_h$-submodule of $E$. Note
that $\G_h$ is an extension of a cyclic group (the image of $\G_h$ under 
$\G@>>>\G/\la a\ra$) by a cyclic group (the intersection $\G_h\cap\la a\ra$). Let 
$H_0=\{h\in H;xc(h)=h\}$. For any $w\in W$ we have 

$\tr(wxc,E)=\sum_{h\in H_0}\tr(wxc,E_h)$. 
\nl
In particular, if $H_0=\em$ then $\tr(wxc,E)=0$ so that the lemma is clear in this case.
Thus we may assume that $H_0\ne\em$. Let $h\in H_0$. We can find a simple 
$\QQ[W]$-submodule $E_h^0$ of $E^h$ such that $E^h=\bfU\ot_\QQ E^0_h$ as $\bfU[W]$-modules.
By Lemma 39.2 applied to $W,\G_h,E_h^0$ instead of $W,\G,E$ we see that the
$\QQ[W]$-structure on $E_h^0$ extends to a $\QQ[W\G_h]$-module structure. (The hypotheses 
of that lemma are satisfied since the $\bfU[W]$-module structure on $E_h$ extends to a 
$\bfU$-module structure.) For $\g\in\G_h$ let $t_\g\in\Aut_\QQ(E_h^0)$ be the action of 
$\g$ in this $\QQ[W\G_h]$-module. By extension of scalars, $t_\g$ defines an element 
$\tit_\g\in\Aut_\QQ(E_h)$. For $e\in E_h,w\in W$ we have 
$\tit_\g(we)=(\g w\g\i)(\tit_\g(e))$. Hence $\g\i(\tit_\g(we))=w\g\i\tit_g(e)$. Thus 
$\g\i\tit_\g:E_h@>>>E_h$ commutes with the action of $W$. Since $E_h$ is simple we see that
$\g\i\tit_\g:E_h@>>>E_h$ is a scalar $\l_h(\g)\in\bfU^*$. Thus 
$\tit_\g=\l_h(\g)\g:E_h@>>>E_h$. Clearly $\g\m\l_h(\g)$ is a homomorphism 
$\G_h@>>>\bfU^*$. For $w\in W$ we have $wxc=\l_h(xc)\i w\tit_{xc}:E_h@>>>E_h$ hence 

$\tr(wxc,E_h)=\l_h(xc)\i\tr(w\tit_{xc},E_h)=\l_{xc}\i\tr(wt_{xc},E_h^0)$.
\nl
Note that $\tr(wt_{xc},E_h^0)\in\ZZ$ since $E_h^0$ is a $\QQ[W\G_h]$-module. We deduce

$\tr(wxc,E)=\sum_{h\in H_0}\l_h(xc)\i\tr(wt_{xc},E_h^0)$.
\nl
Let $H_1=\{h\in H_0;\tr(w't_{xc},E_h^0)\ne 0 \text{ for some }w'\in W\}$. Then, clearly,

$\tr(wxc,E)=\sum_{h\in H_1}\l_h(xc)\i\tr(wt_{xc},E_h^0)$.
\nl
Since $\l_h(xc)$ is a root of $1$, it is enough to verify the following statement:

If $h,h'\in H_1$ then $\l_h(xc)=\pm\l_{h'}(xc)$.
\nl
Since $\la a\ra$ acts transitively on $H$, we can find $y\in\la a\ra$ such that $h=yh'$. We
have $yE_{h'}=E_h,y\G_{h'}y\i=\G_h$. Also, for any $w'\in W$ we have 

$\tr(w'xc,E_{h'})=\tr(yw'xcy\i,E_h)=\tr(yw'y\i yxcy\i,E_h)$
\nl
hence

$\l_{h'}(xc)\tr(w'xc,E_{h'}^0)=\l_h(yxcy\i)\tr(yw'y\i w'xcy\i,E^0_h)$.
\nl
We take here $w'\in W$ such that $\tr(w'xc,E_{h'}^0)\ne 0$. Dividing by 
$\tr(w'xc,E_{h'}^0)$ we deduce $\l_{h'}(xc)\in\QQ\l_h(yxcy\i)$. Since 
$\l_{h'}(xc),\l_h(yxcy\i)$ are roots of $1$, it follows that $\l_{h'}(xc)=\pm\l_h(yxcy\i)$.
Thus it is enough to show that $\l_h(xc)=\pm\l_h(yxcy\i)$ or that 
$\l_h(c\i x\i yxcy\i)=\pm 1$ or that $\l_h(c\i ycy\i)=\pm 1$. More generally, we will 
verify the following statement:

If $h\in H_0,y\in\la a\ra$, $c\i ycy\i\in\la a\ra_h$ then $\l_h(c\i ycy\i)=\pm 1$.
\nl
If $cy=yc=1$ this is obvious. Therefore we may assume that $cy'=y'{}\i c$ for any 
$y'\in\la a\ra$. 

Let $u\in\la a\ra_h$. Since $u,xc$ belong to $\G_h$ we have $\l_h((xc)\i u xc u\i)=1$ that
is, $\l_h(c\i uc u\i)=1$ (we use $ux=xu$). But $c\i uc=u\i$. Hence $\l_h(u^{-2})=1$ that 
is, $\l_h(u)^{-2}=1$ and $\l_h(u)=\pm 1$. We apply this with 
$u=c\i y_1cy_1\i\in\la a\ra_h$. We see that $\l_h(c\i y_1cy_1\i)=\pm 1$. The lemma is 
proved.

\subhead 39.6\endsubhead
Let $W,I,\cx$ are as in 39.1. We assume that $(W,I)$ is irreducible. Let 
$\tI=I\sqc\{\o\}$ where $\o$ is a symbol. We define a map $\p:\tI@>>>W$ as follows: 
$\p(s)=s$ if $s\in I$ and $\p(\o)$ is the unique reflection of maximal length in $\cx$. The
restriction of $\p$ to $\tI$ is injective if $|I|\ge 2$; if $|I|=1$, it maps both elements
of $\tI$ to the unique element of $I$. Let $\Om$ be the group of all permutations 
$\s:\tI@>\si>>\tI$ such that there exist $w\in W$ with $w\p(\s)w\i=\p(\s(x))$ for all 
$x\in\tI$. Then 

$\s\m[w_1\m ww_1w\i]$
\nl  
is a homomorphism of $\Om$ into $\text{\rm Inn}(W)$, the group of inner automorphisms of 
$W$. Let $K$ be a subset of $\tI$ such that $K\ne\tI$. Then $\p$ restricts to a bijection 
$K@>\si>>\p(K)$. Let $\Om^K=\{\s\in\Om;\s(K)=K\}$. Let $W^{(K)}$ be the subgroup of $W$ 
generated by $\p(K)$. From the theory of affine Weyl groups we see that $W^{(K)}$ is a 
(finite) Coxeter group on the generators $\p(K)$. We have canonically 
$W^{(K)}=\prod_{z\in Z}W^{(K)}_z$  where $W^{(K)}_z$ is an irreducible Weyl group with set
of simple reflections $K_z=K\cap W^{(K)}_z$. For $z\in Z$ we set $\cx^K_z=\Rf(W^{(K)}_z)$ 
if $\Rf(W^{(K)}_z)$ is a single $W^{(K)}_z$-conjugacy class and $\cx^K_z=\cx\cap W^{(K)}_z$
if $\Rf(W^{(K)}_z)$ is a union of two $W^{(K)}_z$-conjugacy classes; in any case, $\cx^K_z$
is a single $W^{(K)}_z$-conjugacy class in $\Rf(W^{(K)}_z)$. (We use the following fact: if
$s\in I$ is such that $s\p(\o)$ has order $\ge 4$ then $s,\p(\o)$ are not conjugate under 
$W$.) Then $\cx^K=\sqc_{z\in Z}\cx^K_z$ is a special subset of $\Rf(W^{(K)})$ (see 39.1). 
Now the image of $\Om^K$ under $\Om@>>>\text{\rm Inn}(W)$ (as above) is contained in the 
group $\Aut(W^{(K)},K,\cx^K)$ of automorphisms of $W^{(K)}$ which preserve $K$ and $\cx^K$.
Thus we have a homomorphism $\Om^K@>>>\Aut(W^{(K)},K,\cx^K)$. Restricting this to a 
subgroup $C$ of $\Om^K$ we obtain a homomorphism $C@>>>\Aut(W^{(K)},K,\cx^K)$. Let 
$c\in\Aut(W,I,\cx)$ be such that $c(K)=K$. We have $c^N=1$ for some $N\ge 1$. We extend the
bijection $c:I@>\si>>I$ (restriction of $c:W@>>>W$) to a bijection $c:\tI@>\si>>\tI$ by 
$c(\o)=\o$. We have $c(W^{(K)})=W^{(K)},c(\cx^K)=\cx^K$ hence $c$ restricts to an element 
of $\Aut(W^{(K)},K,\cx^K)$. For any $\s\in\Om$ we define an element $c(\s)\in\Om$ by 
$c(\s)(x)=\s(c\i x)$ for $x\in\tI$. Then $c:\Om@>\si>>\Om$ preserves $\Om^K$. Assume that
$c(C)=C$. For $\s\in C,w\in W^{(K)}$ we have $c(\s(w))=c(\s)(c(w))$. Let $\la c\ra$ be a 
cyclic group of order $N$ with generator $c$. On the set $W^{(K)}\T C\T\la c\ra$ we have a
group structure

$(w,\s,c^n)(w',\s',c^{n'})=(w\s(c^n(w')),\s c^n(\s'),c^{n+n'})$.
\nl
This is the semidirect product $W^{(K)}\G$ of $W^{(K)}$ with the group $\G=C\T\la c\ra$ 
with group structure $(\s,c^n)(\s',c^{n'})=(\s c^n(\s'),c^{n+n'})$. 

\proclaim{Proposition 39.7}In the setup of 39.6 let $E$ be a simple 
$\bfU[W^{(K)}\G]$-module. Let $\g\in\G$. There exists $\z$, a root of $1$ in $\bfU$, such 
that $\tr(w\g,E)\in\z\ZZ$ for any $w\in W^{(K)}$.
\endproclaim
From the theory of affine Weyl groups it is known that one of the following holds:

(i) $\Om$ is cyclic and $c(\s)=\s$ for all $\s\in\Om$.

(ii) $\Om$ is cyclic and $c(\s)=\s\i$ for all $\s\in\Om$,

(iii) $\Om\cong\ZZ/2\T\ZZ/2$.
\nl
Moreover, $C\sub\Om$ (compatibly with the action of $c$). In cases (i),(ii), $C$ is cyclic
and the assumptions of 39.5 are satisfied (with $W^{(K)}$ instead of $W$). Hence the result
follows from 39.5. 

In the remainder of the proof we assume that we are in case (iii). If $C$ is cyclic then
$c(\s)=\s=\s\i$ for all $\s\in C$. The result follows again from 39.5. Hence we may assume
that $C$ is not cyclic so that $C=\Om$. Let $n$ be the order of $c:C@>>>C$. Then 
$n\in\{1,2,3\}$ and $n$ divides $N$. Note that $c^{N/n}$ is in the centre of $\G$. 
Tensoring $E$ by a suitable one dimensional representation (which is trivial on $W^{(K)}C$)
we may assume that $c^{N/n}$ acts trivially on $E$ hence $E$ factors through the quotient 
of $W^{(K)}\G$ by the subgroup generated by $c^{N/n}$. Hence we  may assume that $N=n$. If
$n=1$ then $\G\cong\ZZ/2\T\ZZ/2$ is a group as in 39.5 and the result follows from 39.5. If
$n=2$ then $\G$ is a dihedral group of order $8$ which is again a group as in 39.5 (an 
extension of $\ZZ/2$ by $\ZZ/4$) and the result follows from 39.5. If $n=3$ then $W$ must 
be of type $D_4$, $W^{(K)}$ is an elementary abelian $2$-group, $E$ is one dimensional and
there exists a homomorphism $\mu:W^{(K)}\G@>>>\bfU^*$ such that $\tr(w\g',E)=\mu(w\g')$ for
all $w\in W^{(K)},\g'\in\G$. Then $\z=\mu(\g)$ is a root of $1$. For any $w\in W^{(K)}$ we 
have $\mu(w\g)=\z\mu(w)$ and $\mu(w)=\pm1$, since $w^2=1$. Thus, $\tr(w\g,E)=\pm\z$. The 
proposition is proved.

\proclaim{Corollary 39.8}Assume that $G,D,A,u$ are as in 35.22 and that $G^0/\cz_{G^0}$ is
simple. Then $b_{A,u}^v\in\et\QQ$ for some $\et$, a root of $1$ (with $b_{A,u}^v$ as in 
34.19.)
\endproclaim
This follows from 35.22 and 39.7.


\widestnumber\key{BBD}
\Refs  
\ref\key{\AL}\by D. Alvis\paper The duality operation in the character ring of a finite 
Chevalley group\jour Bull. Amer. Math. Soc.\vol1\yr1979\pages907-911\endref
\ref\key{\BBD}\by A. Beilinson, J. Bernstein, P. Deligne\paper Faisceaux pervers\jour 
Ast\'erisque\vol100\yr1982\endref
\ref\key{\CU}\by C.W. Curtis\paper Truncation and duality in the character ring of a finite
group of Lie type\jour J. Algebra\vol62\yr1980\pages320-332\endref
\ref\key{\DE}\by P. Deligne\paper La conjecture de Weil,II\jour Publ.Math. I.H.E.S.\vol52
\yr1980\pages 137-252\endref
\ref\key{\KA}\by N. Kawanaka\paper Fourier transforms of nilpotently supported invariant 
functions on a simple Lie algebra over a finite field\jour Invent.Math.\vol69\yr1982\pages
411-435\endref
\ref\key{\CRG}\by G. Lusztig\book Characters of reductive groups over a finite field, Ann.
Math. Studies 107\publ Princeton U.Press\yr1984\endref
\ref\key{\CS}\by G. Lusztig\paper Character sheaves,I\jour Adv. Math.\vol56\yr1985\pages
193-237\moreref II\vol57\yr1985\pages226-265\moreref III\vol57\yr1985\pages266-315
\moreref IV\vol59\yr1986\pages1-63\moreref V\vol61\yr1986\pages103-155\endref
\ref\key{\AD}\by G. Lusztig\paper Character sheaves on disconnected groups,I\jour 
Represent. Th. (electronic)\vol7\yr2003\pages374-403\moreref II\vol8\yr2004\pages72-124
\moreref III\vol8\yr2004\pages125-144\moreref IV\vol8\yr2004\pages145-178\moreref Errata
\vol8\yr2004\pages179-179\moreref V\vol8\yr2004\pages346-376\moreref VI\vol8\yr2004\pages
377-413\moreref VII\vol9\yr2005\pages209-266\endref
\ref\key{\PCS}\by G. Lusztig\paper Parabolic character sheaves,I\jour Moscow Math.J.\vol4
\yr2004\pages153-179\endref
\endRefs
\enddocument